%AMSPPT file with "Geometry And Topology Of Complex Hyperbolic And 
%CR-Manifolds "
\magnification=1200
\input amstex
\input psfig
\documentstyle{amsppt}
\NoBlackBoxes
\CenteredTagsOnSplits
\hcorrection{.25in}
\advance\voffset.25in
\advance\vsize-.5in

\TagsOnRight
%\leftheadtext{Boris Apanasov}
\rightheadtext{Complex Hyperbolic and CR-manifolds}
\define\Ga{\Gamma}
\define\ga{\gamma}
\define\a{\alpha}
\define\da{\delta}

\define\e{\epsilon}
\define\La{\Lambda}

\define\sa{\sigma}
\define\Sa{\Sigma}
\define\Om{\Omega}

\define\z{\zeta}
\define\vp{\varphi}
\define\vP{\varPhi}

\define\gtS{{\frak S}}
\define\gts{{\frak s}}

\define\sch{{\Cal H}}

\define\sci{{\Cal I}}

\define\sct{{\Cal T}}
\define\ba{{\Bbb A}}
\define\bb{{\Bbb B}}
\define\bc{{\Bbb C}}

\define\bh{{\Bbb H}}

\define\bp{{\Bbb P}}

\define\br{{\Bbb R}}

\define\bz{{\Bbb Z}}

\define\ch{\bh_{\bc}^}
\define\rh{\bh_{\br}^}
\define\cb{\bb_{\bc}^n}

\define\im{\operatorname{Im}}

\define\is{\operatorname{Isom}}
\define\fix{\operatorname{fix}}

\define\Hom{\operatorname{Hom}}
\define\Vol{\operatorname{Vol}}

\define\ins{\operatorname{int}}

\define\p{\partial}
\define\bs{\backslash}

\define\ra{\rightarrow}
\define\hra{\hookrightarrow}
\define\col{\,:\,}

\define\ov{\overline}

 \topmatter

\title Geometry And Topology Of Complex Hyperbolic And 
CR-Manifolds\endtitle
\author Boris Apanasov
\endauthor
\thanks{Research in MSRI was supported in part 
by NSF grant DMS-9022140.}\endthanks
\address Department of Mathematics, University of Oklahoma,
Norman, OK  73019
\endaddress
\email Apanasov\@ou.edu \endemail
\address 
Sobolev Inst. of Mathematics, Russian Acad. Sci.,
Novosibirsk, Russia 630090 \hfill \hfill \hfill \endaddress
\address 
 Mathematical Sciences Research Institute, 
Berkeley, CA 94720-5070 \hfill \hfill \hfill \endaddress

\subjclass 57, 55, 53, 51, 32, 22, 20 \endsubjclass

\abstract\nofrills {ABSTRACT.} We study geometry, topology and 
deformation spaces of noncompact complex hyperbolic manifolds (geometrically finite, 
with variable negative curvature), whose properties make them surprisingly different from real
hyperbolic manifolds with constant negative curvature. This study uses an interaction between
K\"ahler geometry of the complex hyperbolic 
space and the contact structure at its infinity (the one-point compactification
of the Heisenberg group), in particular an established 
structural theorem for discrete group actions on nilpotent Lie groups.
%(in the case of the Heisenberg group $\sch_n$).
\endabstract
\endtopmatter
\bigskip

\document

\head 1. Introduction \endhead

This paper presents recent progress in studying topology 
and geometry of complex hyperbolic manifolds $M$ with {\it variable} negative curvature and 
spherical Cauchy-Riemannian manifolds with Carnot-Caratheodory
structure at infinity $M_{\infty}$.  

Among negatively curved manifolds, the class of complex hyperbolic manifolds
occupies a distinguished niche due to several reasons. 
First, such manifolds furnish the simplest examples of negatively curved K\"ahler manifolds,
and due to their complex analytic nature, a broad spectrum of 
techniques can contribute to the study. 
Simultaneously, the infinity of
such manifolds, that is the spherical Cauchy-Riemannian manifolds furnish the 
simplest examples of manifolds with contact structures. 
Second, such manifolds provide simplest examples of negatively curved manifolds not having 
{\it constant} sectional curvature, and already obtained results show 
surprising differences between geometry and topology of such manifolds and corresponding properties
of (real hyperbolic) manifolds with constant negative curvature, see \cite{BS, BuM, EMM, Go1, GM, 
Min, Yu1}. 
Third, such manifolds occupy a remarkable place among rank-one symmetric spaces
in the sense of their deformations: they enjoy the flexibility of low dimensional 
real hyperbolic manifolds (see \cite{Th, A1, A2} and \S7) as well as the rigidity of quaternionic/octionic 
hyperbolic manifolds and higher-rank locally symmetric spaces \cite{MG1, Co2, P}.
Finally, since its inception, the theory of smooth 4-manifolds has relied upon complex surface
theory to provide its basic examples. Nowadays it pays back, and
one can study complex analytic 2-manifolds by using  
Seiberg-Witten invariants, decomposition of 4-manifolds along
homology 3-spheres, Floer homology and new (homology) cobordism
invariants, see \cite{W, LB, BE, FS, S, A9} and \S5. 

Complex hyperbolic geometry is the geometry of the unit ball $\cb$ in
$\bc^n$ with the K\"ahler structure given by the Bergman metric (compare
\cite{CG, Go3}, whose automorphisms are biholomorphic automorphisms
of the ball, i.e., elements of $PU(n,1)$. (We notice that complex hyperbolic 
manifolds with non-elementary fundamental groups are complex hyperbolic in 
the sense of S.Kobayashi \cite{Kob}.) 
Here we study topology and geometry of complex hyperbolic manifolds by using
spherical Cauchy-Riemannian geometry at their infinity.
This CR-geometry is modeled on the one point compactification
of the (nilpotent) Heisenberg group, which appears as
the sphere at infinity of the complex hyperbolic space $\ch n$.
In particular, our study exploits a structural Theorem 3.1 about actions of
discrete groups on nilpotent Lie groups (in particular on the 
Heisenberg group $\sch_n$), which generalizes a Bieberbach theorem for Euclidean spaces
\cite{Wo} and strengthens a result by L.Auslander \cite {Au}. 

Our main assumption on a complex hyperbolic $n$-manifold $M$ is
the geometrical finiteness condition on its fundamental group $\pi_1(M)=G\subset PU(n,1)$, 
which in particular implies that $G$ is finitely generated \cite{Bow} and even
finitely presented, see Corollary 4.5. 
The original definition of a geometrically finite manifold $M$ (due to
L.Ahlfors \cite { Ah}) came from an assumption that $M$ may be decomposed 
into a cell by cutting along a finite number of its totally geodesic
hypersurfaces. The notion of
geometrical finiteness has been essentially used in the case of real hyperbolic 
manifolds (of constant sectional curvature), where
geometric analysis and ideas of Thurston have provided powerful tools for
understanding of their structure, see \cite{BM, MA, Th, A1, A3}. Some of those ideas 
also work in spaces with pinched negative curvature, see \cite{Bow}. 
However, geometric methods based on consideration of finite sided fundamental 
polyhedra cannot be used in spaces of variable curvature, see \S4, and we base our 
geometric description of geometrically finite complex hyperbolic manifolds
on a geometric analysis of their ``thin" ends. This analysis is based on establishing 
a fiber bundle structure on Heisenberg (in general, non-compact) manifolds
which remind Gromov's almost flat (compact) manifolds, see \cite{Gr1, BK}.
 
 As an application of our results on geometrical finiteness, we are able to find finite 
coverings of an arbitrary geometrically finite complex hyperbolic 
manifold such that their parabolic ends have the simplest possible structure, i.e., ends with
either Abelian or 2-step nilpotent holonomy (Theorem 4.9). In another such 
an application, we study an 
interplay between topology and K\"ahler geometry of complex hyperbolic
$n$-manifolds, and topology and Cauchy-Riemannian geometry
of their boundary $(2n-1)$-manifolds at infinity, see
our homology cobordism Theorem 5.4. In that respect, the problem of
geometrical finiteness is very different in complex dimension two, where
it is quite possible that complex surfaces with finitely generated
fundamental groups and ``big" ends at infinity are in fact geometrically
finite. We also note that such non-compact geometrically finite complex hyperbolic 
surfaces have infinitely many smooth structures, see \cite{BE}.

The homology cobordism Theorem 5.4 is also an attempt to
control the boundary components at infinity of complex hyperbolic manifolds.
Here the situation is absolutely different from the real hyperbolic one.
In fact, due to Kohn--Rossi analytic extension theorem in the compact case 
\cite{EMM} and to D.Burns theorem in the case when only one boundary component at
infinity is compact (see also \cite{NR1, Th.4.4}, \cite{NR2}), the whole boundary at
infinity of a complex hyperbolic manifold $M$ of infinite volume is
connected (and the manifold itself is geometrically finite if $\dim_{\bc}M\geq 3$) 
if one of the above compactness conditions holds.
However, if boundary components of $M$ are non-compact, the 
boundary $\p_{\infty} M$ may have arbitrarily many components due to our 
construction in Theorems 5.2 and 5.3. 

The results on geometrical finiteness are naturally linked with the Sullivan's
stability of discrete representations of $\pi_1(M)$ into $PU(n,1)$,
deformations of complex hyperbolic manifolds and 
Cauchy-Riemannian manifolds at their infinity, and equivariant (quasiconformal or
quasisymmetric) homeomorphisms inducing such deformations and isomorphisms of discrete 
subgroups of $PU(n,1)$. Results in these directions are discussed in the last two 
sections of the paper.

First of all, complex hyperbolic and CR-structures are very
interesting due to properties of their deformations, 
rigidity versus flexibility. Namely, finite volume complex hyperbolic 
manifolds are rigid due to Mostow's rigidity \cite{Mo1} (for all locally 
symmetric spaces of rank one). Nevertheless their {\it constant curvature}
analogue, real hyperbolic manifolds are flexible in low dimensions and 
in the sense of quasi-Fuchsian deformations (see our discussion in \S7). 
Contrasting to such a flexibility, complex hyperbolic manifolds share 
the super-rigidity of quaternionic/octionic hyperbolic manifolds (see 
Pansu's \cite{P} and Corlette's \cite{Co1-2} rigidity theorems,
analogous to Margulis's \cite{MG1} super-rigidity in higher rank).
 Namely, due to Goldman's \cite{Go1} local rigidity theorem in dimension 
$n=2$ and its extension \cite{GM} for $n\geq 3$, every nearby discrete 
representation $\rho\col G\to PU(n,1)$ of a cocompact lattice 
$G\subset PU(n-1,1)$ stabilizes a complex totally geodesic subspace 
$\ch {n-1}$ in $\ch n$, and for $n\geq 3$, this rigidity is even global 
due to a celebrated Yue's theorem \cite{Yu1}.
  
One of our goals here is to show that, in contrast to that rigidity of 
complex hyperbolic non-Stein manifolds, complex hyperbolic Stein manifolds 
are not rigid in general. Such a flexibility has two aspects. 
Firstly, we point out that the rigidity condition that the group 
$G\subset PU(n,1)$ preserves a complex
totally geodesic hyperspace in $\ch n$ is essential for local rigidity
of deformations only for co-compact lattices $G\subset PU(n-1,1)$. This
is due to the following our result \cite{ACG}:

\proclaim{Theorem 7.1} Let $G\subset PU(1,1)$ be a {\it co-finite} free lattice whose action
in $\ch 2$ is generated by four real involutions (with fixed mutually tangent 
real circles at infinity).
Then there is a continuous family $\{f_\a \}$, $-\e<\a<\e$, of 
$G$-equivariant homeomorphisms 
in $\ov{\ch 2}$ which induce non-trivial quasi-Fuchsian (discrete
faithful)
representations $f^*_\a \col G\to PU(2,1)$. Moreover, for each $\a\neq 0$, 
any $G$-equivariant homeomorphism of $\ov{\ch 2}$ that induces the 
representation $f^*_\a$ cannot be quasiconformal.
\endproclaim

This also shows the impossibility to extend the Sullivan's quasiconformal stability theorem
\cite{Su2} to that situation, as well as provides the first continuous 
(topological) deformation of a co-finite Fuchsian group $G\subset PU(1,1)$ 
into quasi-Fuchsian groups $G_\a=f_\a G f^{-1}_\a\subset PU(2,1)$ with 
the (arbitrarily close to one) Hausdorff dimension 
$\dim_H\La(G_\a)>1$ of the limit set $\La(G_\a)$, $\a\neq1$, 
compare \cite{Co1}.

Secondly, we point out that the noncompactness condition in our non-rigidity
theorem is not essential, either. Namely, complex hyperbolic 
Stein manifolds homotopy equivalent to their closed totally 
{\it real} geodesic surfaces are not rigid, too. Namely, in complex 
dimension $n=2$, we provide 
a canonical construction of continuous 
quasi-Fuchsian deformations of complex surfaces
fibered over closed Riemannian surfaces, which we call ``complex bendings"
along simple close geodesics. This is the first such
deformations (moreover, quasiconformally induced ones) of complex 
analytic fibrations over a compact base:  
 
\proclaim{Theorem 7.2} Let $G\subset PO(2,1)\subset PU(2,1)$ be a given
(non-elementary) discrete group. Then, for any simple closed geodesic 
$\a$ in the Riemann 2-surface $S=H^2_{\br}/G$ and a sufficiently small $\eta_0>0$, 
there is a holomorphic family of $G$-equivariant quasiconformal homeomorphisms 
$F_{\eta}: \ov{\ch 2} \ra \ov{\ch 2}$, $-\eta_0<\eta<\eta_0$, which
defines the bending (quasi-Fuchsian) deformation 
$\Cal B_\a\col (-\eta_0,\,\eta_0)\ra \Cal R_0(G)$ of the group $G$ along the geodesic 
$\a$, $\Cal B_\a(\eta)=F^*_{\eta}$. 
\endproclaim

The constructed deformations depend on many parameters described 
by the Teichm\"uller space 
$\Cal T(M)$ of isotopy classes of complex hyperbolic structures on $M$, or
equivalently by the Teichm\"uller space $\Cal T(G)=\Cal R_0(G)/PU(n,1)$
of conjugacy classes of discrete faithful representations 
$\rho\in\Cal R_0(G)\subset \Hom (G, PU(n,1))$ of $G=\pi_1(M)$:

\proclaim{Corollary 7.3} Let $S_p=\rh 2/G$ be a closed totally real
geodesic surface of genus $p>1$ in a given complex hyperbolic surface
$M=\ch 2/G$, $G\subset PO(2,1)\subset PU(2,1)$. Then there is an
embedding $\pi\circ\Cal B\col B^{3p-3}\hra \Cal T(M)$ of a real $(3p-3)$-ball
into the Teichm\"uller space of $M$,
defined by bending deformations along disjoint closed geodesics in $M$  and
by the projection $\pi\col \Cal D(M)\ra \Cal T(M)=\Cal D(M)/PU(2,1)$ in
the development space $\Cal D(M)$.
\endproclaim

As an application of the constructed deformations, we
answer a well known question about cusp groups on the boundary of
the Teichm\"uller space $\sct (M)$ of a (Stein) complex hyperbolic surface $M$
fibering over a compact Riemann surface of genus $p>1$ \cite{AG}: 

\proclaim {Corollary 7.12} Let $G\subset PO(2,1)\subset PU(2,1)$ be a
uniform lattice isomorphic to the fundamental group of a closed surface
$S_p$ of genus $p\geq 2$. Then there is a continuous deformation 
$R\col \br^{3p-3}\to \Cal T(G)$ (induced by $G$-equivariant 
quasiconformal homeomorphisms of $\ch 2$) whose boundary group
$G_{\infty}=R(\infty)(G)$ has $(3p-3)$ non-conjugate accidental parabolic
subgroups.
\endproclaim
 
Naturally, all constructed topological deformations are in particular 
geometric realizations of the corresponding (type preserving) 
discrete group isomorphisms, see Problem 6.1. However, as Example 6.7 shows,
not all such type preserving isomorphisms are so good. Nevertheless,
as the first step in solving the geometrization Problem 6.1, 
we prove the following geometric realization theorem \cite{A7}:

\proclaim{Theorem 6.2} Let $\phi : G\ra H$ be a type preserving isomorphism
 of two non-ele\-men\-ta\-ry geometrically finite groups $G,H\subset PU(n,1)$. 
Then there exists a unique equivariant ho\-meo\-mor\-phism $f_{\phi}\col \La(G)\ra \La(H)$ of
their limit sets that induces the isomorphism $\phi$. Moreover, 
if $\La(G)=\p \ch n$, the homeomorphism $f_{\phi}$ is the restriction  of a
hyperbolic isometry $h\in PU(n,1)$.
 \endproclaim

We note that, in contrast to Tukia \cite{Tu} isomorphism theorem in 
the real hyperbolic geometry, one might suspect that in general 
the homeomorphism $f_{\phi}$ has no good metric properties, compare 
Theorem 7.1. This is still one of open problems in complex hyperbolic geometry
(see \S6 for discussions).

\head 2. Complex hyperbolic and Heisenberg manifolds  \endhead
 
We recall some facts
concerning the link between nilpotent geometry of the Heisenberg group, 
the Cauchy-Riemannian geometry (and contact structure) 
of its one-point compactification,
and the K\"ahler geometry of the complex hyperbolic space (compare 
\cite{GP1, Go3, KR}).  

One can realize the complex hyperbolic geometry in the complex projective 
space, 
$$
\ch n =\{[z] \in \bc\bp^n\col \langle z,z\rangle<0\,,\,z\in 
\bc ^{n,1} \}\,,
$$
as the set of negative lines in the Hermitian vector space $\bc ^{n,1}$,
with Hermitian structure given by the indefinite $(n,1)$-form
$\langle z,w\rangle = z_1\overline w_1+\cdots+z_n\overline w_n-z_{n+1}\overline w_{n+1}$.
Its boundary $\p\ch n = \{[z]\in \bc\bp ^{n,1}\col \langle z,z\rangle=0\}$ 
consists of all null lines in $\bc\bp^n$ and is homeomorphic
to  the (2n-1)-sphere $S^{2n-1}$. 

The full group $\is\ch n$ of isometries of $\ch n $ is generated by 
the group of holomorphic automorphisms (= the projective unitary group
$PU(n,1)$ defined by the group $U(n,1)$ of unitary automorphisms of $\bc ^{n,1}$)
together with the antiholomorphic automorphism of $\ch n $ 
defined by the $\bc$-antilinear unitary automorphism of $\bc ^{n,1}$ given by 
complex conjugation $z\mapsto \bar z$. The group $PU(n,1)$ can be embedded in a
linear group due to A.Borel \cite{Bor} (cf. \cite{AX1, L.2.1}), hence
any finitely generated group $G\subset PU(n,1)$
is residually finite and has a finite index torsion free subgroup.
Elements $g\in PU(n,1)$ are of the following three types. 
If $g$ fixes a point in $\ch n$, it is called {\it elliptic}.
If $g$ has exactly one fixed point, and it lies in
$\partial \ch n$, $g$ is called {\it parabolic}. If $g$
has exactly two fixed points, and they lie in $\partial \ch n$,
$g$ is called {\it loxodromic}. These three types exhaust all the possibilities.

There are two common models of complex hyperbolic space $\ch n $ as domains in $\bc^n$,
the unit ball $\Bbb B^n_{\bc}$ and the Siegel domain $\frak S_n$. 
They arise from two affine patches in projective space related to $\ch n $ and its boundary.
Namely, embedding $\bc^n$ onto the affine patch of $\bc\bp^{n,1}$
defined by $z_{n+1}\not=0$ (in homogeneous coordinates) as
$A:\bc^n\ra \bc\bp^n$,\, $z \mapsto [(z,1)]$,
we may identify the unit ball $\Bbb B^n_{\bc}(0,1)\subset\bc^n$
with $\ch n = A(\Bbb B^n_{\bc})$. Here the metric in $\bc^n$
is defined by the standard 
Hermitian form $\langle \langle\, ,\,\rangle\rangle$, and  
the induced metric on ${\Bbb B}^n_{\bc}$ is the Bergman metric 
(with constant holomorphic curvature -1) whose 
sectional curvature is between -1 and -1/4.

The Siegel domain model of $\ch n $ arises from the affine patch complimentary 
to a projective hyperplane $H_{\infty}$ which is tangent to $\p\ch n $
at a point $\infty\in\p\ch n$. For example, taking that point $\infty$ as $(0',-1,1)$ with
$0'\in\bc^{n-1}$ and $H_{\infty}=\{[z]\in \bc\bp^n \col z_n + z_{n+1} = 0\}$,
one has the map $\bold S \col \bc^n \ra \bc\bp^n\bs H_{\infty}$ such that
$$
\left( \matrix
z' \\z_n  \endmatrix\right) \longmapsto\left[ \matrix
z' \\ \frac {1}{2} - z_n \\  \frac {1}{2} + z_n \endmatrix\right]
\quad \text{where} \quad z'=\left( \matrix
z_1 \\ \vdots \\ z_{n-1} \endmatrix\right)\in \bc^{n-1}\,.
$$

In the obtained affine coordinates,
the complex hyperbolic space is identified with the {\it Siegel domain}
$$
\frak S_n = \bold S^{-1}(\ch n) =
\{z\in \bc^n \col z_n + \overline z_n > \langle \langle z', z'\rangle \rangle \}\,,
$$
where the Hermitian form is 
$\langle \bold S(z), \bold S(w) \rangle =\langle \langle z',w'
\rangle \rangle - z_n - \overline w_n $. The automorphism group of this 
affine model of $\ch n$  is the group of affine transformations of $\bc^n$ 
preserving $\frak S_n$. Its unipotent radical is the {\it Heisenberg group}
$\sch_n$ consisting of all {\it Heisenberg translations}
$$
T_{\xi,v}\col (w',w_n)\mapsto \left(w'+\xi, w_n+\langle\langle\xi,w'\rangle\rangle+
\frac{1}{2}(\langle\langle\xi,xi\rangle\rangle -iv)\right)\,,
$$
where $w',\xi\in\bc^{n-1}$ and $v\in\br$.

  In particular $\sch_n$ acts simply transitively on $\p\ch n\bs\{\infty\}$, and 
one obtains the {\it upper half space model} for complex hyperbolic space $\ch n$ by 
identifying $ \bc ^{n-1}\times \br\times [0,\infty)$
and $\ov{\ch n}\bs \{\infty\}$ as
$$
(\xi,v,u)\longmapsto  \left[\matrix
\xi \\
\frac {1}{2}(1-\langle\langle \xi,\xi\rangle\rangle -u+iv) \\
\frac {1}{2}(1+\langle\langle \xi,\xi\rangle\rangle +u-iv)
\endmatrix \right]\,,
$$
where $(\xi,v,u)\in \bc ^{n-1}\times \br\times [0,\infty)$ are the horospherical
coordinates of the corresponding
point in $\ov{\ch n}\bs \{\infty\}$ 
(with respect to the point $\infty \in \p \ch n $, see \cite{GP1}).

We notice that, under this identification, the horospheres in $\ch n$ centered at
$\infty$ are the horizontal slices
$H_t=\{(\xi,v,u)\in  \bc ^{n-1}\times \br \times \br_+\col u=t\}$,
and the geodesics running to $\infty$ are the vertical lines
$c_{\xi,v}(t)=(\xi,v,e^{2t})$ passing through points $(\xi,v)\in
\bc ^{n-1}\times \br$. Thus we see that, via the geodesic perspective 
from ${\infty}$, various horospheres correspond as
$ H_t\ra H_u$ with $(\xi,v,t)\mapsto (\xi,v,u)$.

The ``boundary plane"
$\bc ^{n-1}\times \br\times \{0\}=\p\ch n\bs\{\infty\}$ and the horospheres
$H_u=\bc ^{n-1}\times \br\times \{u\}$, $0<u<\infty$, centered at $\infty$ are
identified with the Heisenberg group $\sch_n=\bc^{n-1}\times \br$. 
It is a 2-step
nilpotent group with center $\{0\}\times \br \subset \bc ^{n-1}\times \br$,  
with the isometric action
on itself and on $\ch n$ by left translations:
$$
T_{(\xi_0,v_0)} \col (\xi,v,u)\longmapsto (\xi_0+\xi\,,v_0+v+2\im
\langle\langle \xi_0,\xi\rangle\rangle \,,u)\,,
$$
and the inverse of $(\xi,v)$ is  $(\xi,v)^{-1}=(-\xi,-v)$.
The unitary group $U(n-1)$  acts on  $\sch_n$
and $\ch n $ by rotations:
$A(\xi,v,u)=(A\xi\,,v\,,u)$ for $A\in U(n-1)$. The semidirect product
$\sch(n)=\sch_n \rtimes U(n-1)$ is naturally embedded in
$U(n,1)$ as follows:
$$
A \longmapsto \pmatrix
A&0&0 \\ 0&1&0 \\ 0&0&1 \endpmatrix \in U(n,1)\quad \text{for} \quad A\in U(n-1)\,,
$$
$$
(\xi,v)\longmapsto \pmatrix
I_{n-1}&\xi & \xi\\
-\bar{\xi}^t&1-{1\over 2}(| \xi | ^2-iv)&
-{1\over 2}(|\xi|^2-iv)\\
\bar\xi^t&{1\over 2}(|\xi|^2-iv)&
1+{1\over 2}(|\xi|^2-iv)\endpmatrix \in U(n,1)
$$
where $ (\xi,v)\in \sch_n= \bc ^{n-1}\times \br$ and $\bar \xi^t$ is
the conjugate transpose of $\xi$.

The action of $\sch(n)$ on $\overline{\ch n} \bs \{ \infty \}$
also preserves the Cygan metric $\rho_c$ there, which 
plays the same role as the Euclidean metric does on the upper half-space 
model of the real hyperbolic space $\bh^n$ and is induced by the
following norm:
$$
||(\xi,v,u)||_c=|\,||\xi||^2 + u - iv|^{1/2}\,,\quad 
(\xi,v,u)\in\bc ^{n-1}\times \br\times [0,\infty)\,.\tag2.1
$$

The relevant geometry on each horosphere $H_u\subset \ch n$, 
$H_u\cong\sch_n=\bc ^{n-1}\times \br$, is 
the spherical $CR$-geometry
induced by the complex hyperbolic structure.
 The geodesic perspective from ${\infty}$
defines $CR$-maps between horospheres, which extend to $CR$-maps
between the one-point compactifications 
$H_u\cup {\infty}\approx S^{2n-1}$.
In the limit, the induced  metrics on horospheres fail to converge 
but the $CR$-structure
remains fixed. In this way, the complex hyperbolic geometry induces
$CR$-geometry on the sphere at infinity $\p \ch n \approx S^{2n-1}$, naturally
identified with the one-point compactification of the Heisenberg group
$\sch_n$.

\head 3. Discrete actions on nilpotent groups and Heisenberg manifolds\endhead
 
In order to study the structure of Heisenberg manifolds (i.e., the manifolds 
locally modeled on the Heisenberg group $\sch_n$)
and cusp ends of complex hyperbolic
manifolds, we need a Bieberbach type structural theorem for
isometric discrete group actions on $\sch_n$, originally proved in \cite{AX1}.
It claims that each discrete isometry group of
the Heisenberg group $\sch _n$ preserves some left coset of a connected
Lie subgroup, on which the group action is cocompact. 

Here we consider more general situation. Let $N$ be a connected, simply connected 
nilpotent Lie group, $C$ a 
compact group of automorphisms of $N$, and $\Ga$ a discrete subgroup of the 
semidirect product $N\rtimes C$. Such discrete groups are the holonomy groups 
of parabolic ends of locally symmetric rank one (negatively curved) 
manifolds and can be described as follows.

\proclaim{Theorem 3.1} There exist a connected Lie subgroup
$V$ of $N$ and a finite index normal subgroup $\Gamma^*$ 
of $\Gamma$ with the following properties:
\roster
\item There exists $b\in N$ such that $b\Gamma b^{-1}$ preserves
$V$.
\item $V/b\Gamma b^{-1}$ is compact.
\item $b\Gamma^* b^{-1}$ acts on $V$ by left translations and this action is free.
\endroster
\endproclaim

\remark{Remark 3.2} (1) It immediately follows that any discrete subgroup 
$\Ga\subset N\rtimes C$ is virtually nilpotent because it has a finite index 
subgroup $\Gamma^*\subset \Ga$ isomorphic to a lattice in $V\subset N$.

(2) Here, compactness of $C$ is an essential condition
because of Margulis \cite {MG2} construction of nonabelian free 
discrete subgroups $\Gamma$ of $R^3\rtimes GL(3,R)$.

(3) This theorem generalizes a Bieberbach theorem for Euclidean spaces, see \cite{Wo},
and strengthens a result by L.Auslander \cite {Au} who claimed those 
properties not for whole group $\Ga$ but only for its  finite index subgroup. 
Initially in \cite {AX1}, we proved this theorem for the Heisenberg group 
$\sch_n$ where we used Margulis Lemma \cite {MG1, BGS} and geometry 
of $\sch_n$ in order to extend the classical arguments in \cite {Wo}. 
In the case of general nilpotent groups, our proof uses different ideas
and goes as follows (see\cite{AX2} for details). 
\endremark 

\demo{Sketch of Proof} Let $p:\Gamma\rightarrow C$ be the 
composition of the inclusion $\Gamma\subset N\rtimes C$ and the projection
$N\rtimes C\rightarrow C$, $G$ the identity component of 
$\overline{\Gamma N}$, and $\Gamma_1=G\cap \Gamma$. 
Due to compactness of $C$, $G$ has finite index in 
$\overline {\Gamma N}$, so $\Gamma_1$ has finite index in $\Gamma$. Let
$W\subset N$ be the analytic subgroup pointwise fixed 
by $p(\Gamma_1)$. Due to \cite{Au}, for all 
$\gamma=(w,c)\in\Gamma_1$, $w$ lies in $W$.  Thus $\Gamma$ preserves $W$
and, by replacing $N$ with $W$, we may assume that $\Phi=p(\Gamma)$ is 
finite.

   Consider $\Gamma^*=\ker (p)$ which
is a discrete subgroup of $N$ and has finite index in $\Gamma$.
Let $V$ be the connected Lie subgroup of $N$ in which $\Gamma^*$
is a lattice. Then the conjugation action of $\Gamma$ on 
$\Gamma^*$ induces a $\Gamma$-action on $V$. We form the
semi-direct product
$V\rtimes \Gamma$ and let $K=\{(a^{-1},(a,1))\in V\rtimes \Gamma : 
(a,1)\in\Gamma^*\}$. Obviously, $K$ is a normal subgroup of 
$V\rtimes \Gamma$. Defining the maps
$i:V\rightarrow V\rtimes \Gamma/K$ by $i(v)=(v,(1,1))K$ 
and
$\pi:V\rtimes \Gamma/K\rightarrow \Phi$ by 
$\pi(v,(a,A))=A$, we get a short exact sequence
$$\CD
1  @>>> V  @>i>> V\rtimes \Gamma/K @>\pi>> \Phi @>>>1\,.
\endCD
$$

Since any extension of a finite group by a simply connected nilpotent 
Lie group splits, there is a homomorphism
$s: \Phi\rightarrow V\rtimes \Gamma/K$ such that 
$\pi\circ s=id_\Phi$.
For each $A\in \Phi$, we fix an element $(f(A),(g(A),A))\in V\rtimes\Ga$
representing $s(A)$. Since $s$ is a homomorphism, we have
$$
g(AB)^{-1}f(AB)^{-1}=A\left(g(B)^{-1}f(B)^{-1}\right)g(A)^{-1}f(A)^{-1}
\quad \text{for}\quad A,B\in \Phi\,.\tag4.3
$$

 Define $h: \Phi\rightarrow N$ by $h(A)=g(A)^{-1}f(A)^{-1}$.
Then (2.4) shows that $h$ is a cocycle. 
Since $\Phi$ is finite and $N$ is a simply connected nilpotent Lie group,
$H^1(\Phi,N)=0 $ due to \cite {LR}. Thus there exists $b\in N$ such that 
$h(A)=A(b^{-1})b$ for all $A\in \Phi$.

On the other hand, $\pi((1,(a,A))K)=\pi((f(A),(g(A),A))K)=A$
for any $\gamma=(a,A) \in \Gamma$. It follows 
that there is $v_0\in V$ such that $ a^{-1}v_0=h(A)$. This and (4.3)
imply that $a^{-1}v_0=A(b^{-1})b$, and hence $baA(b^{-1})=bv_0b^{-1}$.

  Now consider the group $b\Gamma b^{-1}$ which acts on 
$bVb^{-1}$.,
%$b\Gamma b^{-1}\times bV b^{-1}\rightarrow bV b^{-1}$. 
For any $\gamma=(a,A)\in \Gamma$, the action of the element 
$b\gamma b^{-1}=(baA(b^{-1}),A)$ on $bVb^{-1}$ is as follows:
$$
((baA(b^{-1}),A),v')\rightarrow 
baA(b^{-1})A(v')(baA(b^{-1}))^{-1}\,.
$$

In particular, $baA(b^{-1})\cdot A(bVb^{-1})\cdot 
(baA(b^{-1}))^{-1}=bVb^{-1}$. Therefore, $A(bVb^{-1})=bVb^{-1}$ 
because of $baA(b^{-1})=bv_0 b^{-1}\in bVb^{-1}$,
and hence $b\gamma b^{-1}$ preserves $bVb^{-1}$.
\line{\hfil\hfil\hfil\qed}
\enddemo

Now we can apply our description of discrete group actions on a nilpotent group 
(Theorem 3.1) to study the structure of Heisenberg manifolds. Such manifolds are 
locally modeled on the $(\sch_n, \sch (n))$-geometry and each of them 
can be represented as the quotient
$\sch_n/G$ under a discrete, free isometric action of its fundamental
group $G$ on $\sch_n$, i.e., the isometric action of 
a torsion free discrete subgroup of
$\sch(n)=\sch_n\rtimes U(n-1)$. Actually, we establish fiber bundle structures on 
all noncompact Heisenberg manifolds: 

\proclaim{Theorem 3.3} Let $\Ga\subset \sch_n\rtimes U(n-1)$ be a torsion-free 
discrete group  acting on the Heisenberg group $\sch_n= \bc ^{n-1}\times \br $ 
with non-compact quotient.
Then the quotient $\sch_n/ \Ga$ has zero Euler characteristic and
is a vector bundle over a  compact manifold. Furthermore,
this compact manifold is finitely covered
by a nil-manifold which is either a torus or the total space of a 
circle bundle over a torus.
\endproclaim

The proof of this claim (see \cite{AX1}) is based on two facts due to Theorem 3.1. 
First, that the discrete holonomy group $\Ga\cong\pi_1(M)$
of any noncompact Heisenberg manifold $M=\sch_n/\Ga$, $\Ga\subset \sch(n)$, has a 
proper $\Ga$-invariant subspace $\sch_\Ga\subset \sch_n$. And second, 
the compact manifold $\sch_{\Ga}/\Ga$ is  finitely covered by
$\sch_{\Ga}/\Ga^*$ where $\Ga^*$ acts on $\sch_{\Ga}$ by translations. The structure of
the covering manifold $\sch_G/G^*$ is given in the following lemma.

\proclaim{Lemma 3.4} Let $V$ be a connected Lie subgroup of the
Heisenberg group $\sch_n$ and $G\subset V$ a discrete co-compact subgroup of $V$.
 Then the manifold $V/G$ is
\roster
\item  a torus if $V$ is Abelian;
\item  the total space of a torus bundle over a torus if $V$ is not Abelian.
\endroster
\endproclaim

Though noncompact Heisenberg manifolds $M$ are
vector bundles $\sch_n/\Ga\rightarrow \sch_\Ga/\Ga$, simple examples 
show \cite{AX1} that such vector bundles may be non-trivial in general.
However, up to finite coverings, they are trivial \cite{AX1}:

\proclaim{Theorem 3.5} Let $\Ga\subset \sch_n\rtimes U(n-1)$ be a
discrete group and $\sch_\Ga\subset \sch_n$  a 
connected $\Ga$-invariant Lie subgroup on which $\Ga$ acts co-compactly.
Then there exists a finite index subgroup $\Ga_0\subset \Ga$
such that the vector bundle $\sch_n/\Ga_0\rightarrow \sch_\Ga/\Ga_0$ is
trivial. In particular, any Heisenberg orbifold $\sch_n/\Ga$ is
finitely covered by the product of
a compact nil-manifold $\sch_\Ga/\Ga_0$ and an Euclidean space.
\endproclaim

We remark that in the case when $\Ga\subset \sch_n\rtimes U(n-1)$ is a lattice, 
that is the quotient $\sch_n/\Ga$ is compact,
the existence of such finite cover of $\sch_n/\Ga$ by a closed nilpotent
manifold $\sch_n/\Ga_0$ is due to Gromov \cite {Gr} and Buser-Karcher \cite{BK}
results for almost flat manifolds. 

Our proof of Theorem 3.5 has the following scheme.
Firstly, passing to a finite index subgroup, 
we may assume that the group $\Ga$
is torsion-free. After that, we shall find a finite index subgroup
$\Ga_0\subset \Ga$ whose rotational part is ``good". Then we shall
 express the vector bundle $\sch_n/\Ga_0\rightarrow \sch_\Ga/\Ga_0$ as the
 Whitney sum of a trivial bundle and a fiber product. We finish the proof
by using the following criterion about the triviality of fiber products:

\proclaim{Lemma 3.6}
Let $F\times_H V$ be a fiber product and suppose that the homomorphism
$\rho\col H\rightarrow GL(V)$ extends to a homomorphism
$\rho\col F\rightarrow GL(V)$. Then $F\times_H V$  is a trivial bundle,
$F\times_H V\cong F/H\times V$.
\endproclaim

\demo{Proof} The isomorphism $F\times_H V\cong F/H\times V$ is given by
$[f,v]\rightarrow (Hf,\rho(f)^{-1}(v))$.
%\line{\hfil \hfil 
\qed %}
\enddemo

\head 4. Geometrical finiteness in complex hyperbolic geometry  \endhead

Our main assumption on a complex hyperbolic $n$-manifold $M$ is
the geometrical
finiteness of its fundamental group $\pi_1(M)=G\subset PU(n,1)$, which in
particular implies that the discrete group $G$ is finitely generated. 

Here a subgroup $G\subset PU(n,1)$ is called {\it discrete} if it is a 
discrete subset of $PU(n,1)$. The {\it limit set} $\La(G)\subset \p \ch n$ 
of a discrete group $G$ is the set of
accumulation points of (any) orbit $G(y),\, y\in \ch n$. The complement of
$\La(G)$ in $\p \ch n$ is called the {\it discontinuity set} $\Om(G)$.
A discrete group $G$ is called {\it elementary}
if its limit set $\La(G)$ consists of at most two points.  
An infinite discrete group $G$ is called {\it parabolic}
if it has exactly one fixed point $\fix (G)$; then $\La(G)=\fix (G)$, and
$G$ consists of either parabolic or elliptic elements. As it was
observed by many authors (cf. \cite {MaG}), parabolicity in the
variable curvature case is not as easy a condition to deal with as it is
in the constant curvature space. However the results of \S 2 
 simplify the situation, especially for geometrically finite groups.

Geometrical finiteness has been essentially used for real hyperbolic manifolds, where
geometric analysis and ideas of Thurston provided powerful tools for
understanding of their structure. Due to the absence of totally geodesic 
hypersurfaces in a space of variable negative curvature,  we cannot use the original
definition of geometrical finiteness which came from
an assumption that the corresponding real hyperbolic manifold $M=\bh^n/G$ may be 
decomposed  into a cell
by cutting along a finite number of its totally geodesic hypersurfaces, that is
the group $G$ should possess a finite-sided fundamental polyhedron, see  
\cite{Ah}.  However, we can define {\it geometrically finite} groups $G\subset PU(n,1)$
as those ones whose limit sets $\La(G)$ consist of only conical limit points and
parabolic (cusp) points $p$ with compact quotients $(\La(G)\bs\{p\})/G_p$
with respect to parabolic stabilizers $G_p\subset G$ of $p$, see
\cite {BM, Bow}. There are other definitions of geometrical finiteness 
in terms of ends and the minimal convex retract of the noncompact manifold $M$, 
which work well not only in the real hyperbolic spaces $\bh^n$ 
(see \cite {Mar, Th, A1, A3}) but also in spaces with variable pinched negative 
curvature \cite {Bow}.

Our study of geometrical finiteness in complex hyperbolic geometry 
is based on analysis of geometry and topology of thin (parabolic) ends of 
corresponding manifolds and parabolic cusps of discrete isometry groups $G\subset PU(n,1)$. 

Namely, suppose a point $p\in \p \ch n$ is fixed by some
parabolic element of a given discrete group $G\subset PU(n,1)$, and 
$G_p$ is the stabilizer of $p$ in $G$. 
Conjugating $G$ by an element $h_p\in PU(n,1)\,, h_p(p)=\infty$,
we may assume that the stabilizer $G_p$ is a subgroup $G_{\infty}\subset \sch(n)$.
In particular, if $p$ is the origin $0\in\sch_n$, the transformation $h_0$ can be 
taken as the Heisenberg inversion $\sci$ in the hyperchain $\p\ch {n-1}$.
It preserves the unit Heisenberg sphere 
$S_c(0,1)=\{(\xi,v)\in\sch_n\col ||(\xi,v)||_c=1\}$ and acts in 
$\sch_n$ as follows:
$$
\sci(\xi,v)=\left (\frac {\xi}{|\xi|^2-iv}\,,\, \frac {-v}{v^2+|\xi|^4}
\right )\quad \text {where}\,\,(\xi,v)\in \sch_n=\bc^{n-1}\times \br \,.
\tag4.1
$$ 

For any other point $p$, we may take $h_p$ as the Heisenberg 
inversion $\sci_p$ which preserves the unit Heisenberg sphere 
$S_c(p,1)=\{(\xi,v)\col \rho_c(p,\,(\xi,v))=1\}$ centered at $p$.
The inversion $\sci_p$ is conjugate of $\sci$ by the Heisenberg
translation $T_p$ and maps $p$ to $\infty$.

 After such a conjugation, we can apply Theorem 3.1 to the
parabolic stabilizer $G_{\infty}\subset \sch(n)$ and get
a connected Lie subgroup $\sch_\infty
\subseteq \sch_n$ preserved by $G_{\infty}$ (up to changing the origin). So we can
make the following definition.

\definition{Definition 4.2} A set  $U_{p,r}\subset
\overline{\ch n}\bs \{p\}$ is called a {\it standard cusp neighborhood
of radius} $r>0$ at a parabolic fixed point $p\in \p \ch n$ of a discrete
group $G\subset PU(n,1)$ if, for the Heisenberg inversion 
$\sci_p\in PU(n,1)$ with
respect to the unit sphere $S_c(p,1)$,
$\sci_p(p)=\infty\,,$ the following conditions hold:
\roster
\item $U_{p,r}=\sci_p^{-1}\left(\{x\in \ch n\cup \sch_n\col
\rho_c(x,\sch_\infty)\geq 1/r \}\right)\,;$
\item  $U_{p,r}$ is precisely invariant
with respect to $G_p\subset G$, that is:
$$
\gamma(U_{p,r})=U_{p,r} \quad  \text{for} \,\,\;\gamma\in G_p\quad \text {and}\quad
g(U_{p,r})\cap U_{p,r}=\emptyset \quad \text {for} \,\,\,g\in G\backslash
G_p\,.
$$
\endroster
A parabolic point $p\in \p \ch n$ of $G\subset PU(n,1)$ is called 
a {\it cusp point} if it has a cusp neighborhood $U_{p,r}$.
\enddefinition

We remark that some parabolic points of a discrete group
$G\subset PU(n,1)$ may not be cusp points, see examples in \S 5.4 of \cite{AX1}. Applying
Theorem 3.1 and \cite{Bow}, we have:

\proclaim {Lemma 4.3} Let $p\in \p \ch n$ be a parabolic fixed point of
a discrete subgroup $G$ in $PU(n,1)$. Then $p$ is a cusp point if and only if
$(\La(G)\bs \{p\})/G_p$ is compact.
\endproclaim

This and finiteness results of Bowditch \cite {B} allow us to use another 
equivalent definitions of geometrical finiteness.
In particular it follows that a discrete subgroup $G$ in $PU(n,1)$ is
{\it geometrically finite} if and only if its quotient space $M(G)=[\ch n\cup\Om (G)]/G$ 
has finitely many ends, and each of them is a cusp end, that is an end whose neighborhood can be
taken (for an appropriate $r>0$) in the form:
$$
U_{p,r}/G_p\approx (S_{p,r}/G_p) \times (0,1]\,, \tag4.4
$$
where 
$$
S_{p,r}=\p_H U_{p,r}=\sci_p^{-1}\left(\{x\in H_ \bc ^n\cup \sch_n
 \col \rho_c(x,\sch_\infty)=1/r\}\right)\,.
$$

Now we see that a geometrically finite manifold
can be decomposed into a compact submanifold and finitely many cusp
submanifolds of the form (4.4). Clearly, each of such cusp ends is homotopy equivalent
to a Heisenberg $(2n-1)$-manifold and moreover, due to Theorem 3.3, 
to a compact $k$-manifold, $k\leq 2n-1$. From the last fact, it follows that
the fundamental group of a Heisenberg manifold
is finitely presented, and we get the following finiteness result:

\proclaim{Corollary 4.5} Geometrically finite groups
$G\subset PU(n,1)$ are finitely presented.
\endproclaim

In the case of variable curvature, it is problematic to use
geometric methods based on consideration of finite sided
fundamental polyhedra, in particular, Dirichlet
polyhedra $D_y(G)$ for $G\subset PU(n,1)$  bounded by bisectors in a complicated way,
see \cite { Mo2, GP1, FG}. In the case of discrete parabolic
groups $G\subset PU(n,1)$, one may expect that the  Dirichlet
polyhedron $D_y(G)$ centered at a point $y$ lying in a $G$-invariant
subspace has finitely many sides. It is true for real hyperbolic spaces
\cite{A1} as well as for cyclic and dihedral parabolic groups in
complex hyperbolic spaces. Namely, due to \cite{ Ph}, 
Dirichlet polyhedra $D_y(G)$ are always two sided for any cyclic group
$G\subset PU(n,1)$ generated by a Heisenberg translation.
Due the main result in \cite{GP1}, this finiteness also holds
for a cyclic ellipto-parabolic group or a dihedral parabolic group 
$G\subset PU(n,1)$ generated by inversions
in asymptotic complex hyperplanes in $\ch n$ if the central point $y$ lies
in a $G$-invariant vertical line or $\br$-plane (for any other center $y$, 
$D_y(G)$ has infinitely many sides). Our technique easily implies 
that this finiteness still holds for generic parabolic cyclic groups 
\cite{AX1}:

\proclaim{Theorem 4.6}
 For any discrete group $G\subset PU(n,1)$
generated by a parabolic ele\-ment, there exists a point  $y_0\in \ch n$
 such that the Dirichlet polyhedron $D_{y_0}(G)$ centered at $y_0$
has two sides.
 \endproclaim
\demo{Proof} Conjugating $G$ and due to Theorem 3.1, we may assume that $G$ preserves
a one dimensional subspace $\sch_{\infty}\subset \sch_n$ as well as 
$\sch_{\infty}\times\br_+\subset \ch n$, where $G$ acts by translations. So we can
take any point $y_0\in\sch_{\infty}\times\br_+$ as the central point of (two-sided)
Dirichlet polyhedron $D_{y_0}(G)$ because its orbit $G(y_0)$ coincides with the orbit
$G'(y_0)$ of a cyclic group generated by the Heisenberg translation induced by $G$.
\line{\hfil \hfil \qed}
\enddemo

However, the behavior of Dirichlet polyhedra for parabolic groups 
$G\subset PU(n,1)$ of rank more than one can be very bad. It is given by 
our construction \cite{AX1}, where we have evaluated intersections of Dirichlet bisectors with
a 2-dimensional slice:

\proclaim{Theorem 4.7}  Let $G\subset PU(2,1)$ be a discrete parabolic
group conjugate to the subgroup $\Ga=\{(m,n)\in \bc\times \br \col  m,n\in \bz\}$ 
of the Heisenberg group
$\sch_2=\bc \times \br$.
Then any Dirichlet polyhedron $D_y(G)$ centered at any point 
$y\in \ch 2$ has infinitely many sides.
\endproclaim

 Despite the above example, the below application of Theorem 3.1
provides a construction of fundamental polyhedra $P(G)\subset\ch n$
for arbitrary discrete parabolic groups $G\subset PU(n,1)$, which are
bounded by finitely many
hypersurfaces (different from Dirichlet bisectors). 
This result may be seen as a base for extension of
Apanasov's construction \cite {A1} of finite sided
pseudo-Dirichlet polyhedra
in $\bh^n$ to the case of the complex hyperbolic space $\ch n$.

\proclaim{Theorem 4.8}
 For any discrete parabolic group $G\subset PU(n,1)$,
there exists a finite-sided
fundamental polyhedron $P(G)\subset \ch n$.
\endproclaim

 \demo{Proof}
After conjugation, we may assume that
$G\subset \sch_n\rtimes U(n-1)$. Let 
$\sch_{\infty}\subset\sch_n= \bc^{n-1}\times \br$ be the connected $G$-invariant 
subgroup given by 
Theorem 3.1.
For a fixed $u_0>0$, we consider the horocycle
$V_{u_0}=\sch_{\infty}\times \{u_0\}\subset \bc^{n-1}\times \br\times \br_+
=\ch n$. For
distinct points $y,y'\in V_{u_0}$,  the bisector
$C(y,y')=\{z\in\ch n : d(z,y)=d(z,y')\}$ intersects $V_{u_0}$ transversally.
Since
$V_{u_0}$ is $G$-invariant, its intersection with a Dirichlet
polyhedron 
$$
D_y(G)=\bigcap_{g\in G\backslash \{id\} } \{w\in \ch n \col d(w,y)<d(w,g(y))\}
$$
centered at a point $y\in V_{u_0}$ is a
fundamental polyhedron for the $G$-action on $V_{u_0}$. 
The polyhedron $D_{y}(G)\cap V_{u_0}$
is compact due to Theorem 3.3, and hence has finitely many sides.
Now, considering 
$G$-equivariant projections \cite{AX1}:
$$
\pi\col \sch_n\ra\sch_{\infty}\,, \quad \pi'\col \ch n=\sch_n\times \br_+\ra
 V_{u_0}\,,\quad
\pi'(x,u)=(\pi(x),u_0)\,,
$$
we get a finite-sided
fundamental polyhedron ${\pi '}^{-1}(D_y(G)\cap V_{u_0})$ for the action 
of $G$ in $\ch n$.

\line{\hfil \hfil \hfil \qed}
\enddemo

Another important application of Theorem 3.1 shows that
cusp ends of a geometrically finite complex hyperbolic orbifolds $M$ have,
up to a finite covering of $M$, a 
very simple structure:

\proclaim{Theorem 4.9} Let $G\subset PU(n,1)$ be a geometrically
finite discrete group.
Then $G$ has a subgroup $G_0$ of finite
index such that every parabolic subgroup of $G_0$ is isomorphic
to a discrete subgroup of  the Heisenberg group $\sch_n=\bc^{n-1}\times \br$ .
In particular, each parabolic subgroup of $G_0$ is free Abelian or
2-step nilpotent.
\endproclaim

The proof of this fact \cite{AX1} is based on the residual finiteness of geometrically 
finite subgroups in $PU(n,1)$ and the following two lemmas.

\proclaim {Lemma 4.10}  Let $G \subset \sch_n\rtimes U(n-1)$ be a discrete
 group and $\sch_G\subset \sch_n$ a minimal $G$-invariant connected
 Lie subgroup (given by Theorem 3.1).
Then $G$ acts on  $\sch_G$ by translations if $G$ is either
 Abelian or 2-step nilpotent.
\endproclaim

\proclaim{Lemma 4.11} Let $G\subset \sch_n\rtimes U(n-1)$ be a torsion free
discrete group, $F$ a finite group and $\phi\col  G\longrightarrow F$ an
epimorphism.
 Then the rotational part of $ \ker (\phi)$
has strictly smaller order than that of $G$ if one of the following
happens:
 \roster
 \item $G$ contains a  finite index Abelian subgroup and $F$ is not Abelian;
\item $G$ contains a finite index 2-step nilpotent subgroup  and $F$ is not
a 2-step nilpotent group.
\endroster
\endproclaim

We remark that the last Lemma generalizes a result of 
C.S.Aravinda and T.Farrell \cite{AF} for
Euclidean crystallographic groups.

 We conclude this section by pointing out that the problem of geometrical
finiteness is very different in complex dimension two.
Namely, it is a well
known fact that any finitely generated
discrete subgroup of $PU(1,1)$ or $PO(2,1)$
is geometrically finite. This and Goldman's \cite { Go1}
local rigidity theorem
for cocompact lattices $G\subset U(1,1)\subset PU(2,1)$ 
allow us to formulate the following conjecture:

\proclaim{Conjecture 4.12} All finitely generated discrete groups
 $G\subset PU(2,1)$ with non-empty discontinuity set 
$\Om(G)\subset \p\ch 2$ are geometrically finite. 
\endproclaim

\head 5. Complex homology cobordisms and 
the boundary at infinity
\endhead

The aim of this section is to study the
topology of   complex analytic ''Kleinian"
manifolds $M(G) = [\ch n \cup \Om(G)]/G$ with geometrically finite holonomy
groups $G\subset PU(n,1)$. The boundary of this manifold, $\p M=\Om(G)/G$,
has a spherical $CR$-structure and, in general, is non-compact. 

We are especially interested in the case of complex analytic surfaces, 
where powerful methods of 4-dimensional topology may be used.
It is still unknown what are suitable cuts of 
4-manifolds, which (conjecturally) split them into geometric 
blocks (alike Jaco-Shalen-Johannson decomposition of 3-manifolds
in Thurston's geometrization program; for a classification of 
4-dimensional geometries, see \cite{F, Wa}). Nevertheless,
studying of complex surfaces suggests that in this case 
one can use integer homology 3-spheres and ``almost flat" 
3-manifolds (with virtually nilpotent
fundamental groups).
Actually, as Sections 3 and 4 show, the latter manifolds appear at the
ends of finite volume complex hyperbolic manifolds. 
As it was shown by C.T.C.Wall \cite{Wa}, 
the assignment of the appropriate 4-geometry (when available) gives
a detailed insight into the intrinsic structure of a complex surface.
To identify complex hyperbolic blocks in such a splitting, one can use
Yau's uniformization theorem \cite{Ya}. It implies that
every smooth complex projective 2-surface $M$ with positive canonical bundle
and satisfying the topological condition that
$\chi(M)=3\operatorname{Signature}(M)$, is a complex hyperbolic manifold.
The necessity of homology sphere decomposition in dimension
four is due to M.Freedman and L.Taylor result  (\cite{ FT}):
\vskip4pt
{\it Let $M$ be a simply connected 4-manifold with intersection
form $q_M$ which decomposes as a direct sum $q_M=q_{M_1}\oplus q_{M_2}$, where
$M_1, M_2$ are smooth manifolds. Then the manifold $M$ can be represented as
a connected sum $M=M_1\#_{\Sigma}M_2$ along a homology sphere $\Sigma$.}
\vskip4pt
Let us present an example of such a splitting, $M = X \#_{\Sigma} Y$, 
of a simply connected complex surface $M$ 
with the intersection form $Q_M$ into smooth manifolds (with boundary) 
$X$ and $Y$, along a $\bz$-homology 3-sphere $\Sa$ 
such that $Q_M = Q_X \oplus Q_Y$.
Here one should mention that though $X$ and $Y$ are no longer closed
manifolds, the intersection forms $Q_X$ and $Q_Y$ are well defined on the
second cohomology and are unimodular due to the condition that
$\Sigma$ is a $\bz$-homology 3-sphere.

\proclaim{Example 5.1} Let $M$ be the Kummer surface 
$$
K3=\{[z_0,z_1,z_2,z_3]\in \bc\bp^3\col z_0^4+z_1^4+z_2^4+z_3^4=0\}\,.
$$
Then there are four disjointly embedded (Seifert fibered) $\bz$-homology 
3-spheres in $M$, which split the Kummer surface into five blocks: 
$$
K3 = X_1 \cup_{\Sigma} Y_1 \cup_{\Sigma'} Y_2 \cup_{-\Sigma'} Y_3 
\cup_{-\Sigma} X_2\,,
$$
with intersection forms $Q_{X_j}$ and $Q_{Y_i}$ equal $E_8$ and $H$,
respectively:
%and those ones of $Y_i$'s equal $H,\ i = 1,\ldots,3$\,:
$$
E_8= \pmatrix
-2 &  1 &  0 &  0 &  0 &  0 &  0 &  0 \\
 1 & -2 &  1 &  0 &  0 &  0 &  0 &  0 \\
 0 &  1 & -2 &  1 &  0 &  0 &  0 &  0 \\
 0 &  0 &  1 & -2 &  1 &  0 &  0 &  0 \\
 0 &  0 &  0 &  1 & -2 &  1 &  0 &  1 \\
 0 &  0 &  0 &  0 &  1 & -2 &  1 &  0 \\
 0 &  0 &  0 &  0 &  0 &  1 & -2 &  0 \\
 0 &  0 &  0 &  0 &  1 &  0 &  0 & -2
\endpmatrix\,,\quad H = \pmatrix 0&1\\1&0\endpmatrix\,.
$$
Here the $\bz$-homology spheres $\Sigma$ and $\Sigma'$ are correspondingly 
the Poincar\'e homology sphere $\Sigma(2,3,5)$ and
Seifert fibered homology sphere $\Sigma(2,3,7)$; the minus sign means the
change of orientation.
\endproclaim

\demo{Scheme of splitting} Due to J.Milnor \cite{Mil} (see also \cite{RV}), 
all Seifert fibered homology 3-spheres $\Sa$ can be seen as the boundaries at infinity
of (geometrically finite) complex hyperbolic orbifolds $\ch 2/\Ga$, where
the fundamental groups $\pi_1(\Sa)=\Ga\subset PU(2,1)$ act free 
in the sphere at infinity $\p\ch 2=\ov{\sch_2}$. In particular,
the Seifert fibered homology sphere $\Sa'=\Sigma (2,3,7)$ is diffeomorphic to 
the quotient $[(\bc\times\br)\bs(\{0\}\times\br]/\Ga(2,3,7)$. Here
$[(\bc\times\br)\bs(\{0\}\times\br)]$ is the complement in the 3-sphere
$\ov{\sch_2}=\p B^2_{\bc}$ to the boundary circle at infinity of the complex 
geodesic $B^2_{\bc}\cap (\bc\times\{0\})$, and the group
$\Ga(2,3,7)\subset PU(2,1)$ acts on this complex geodesic as
the standard triangle group $(2,3,7)$ in the disk Poincar\'e model 
of the hyperbolic 2-plane $\bh^2_{\br}$. 

This homology 3-sphere $\Sa'$ 
embeds in the $K3$-surface $M$, splitting it into submanifolds with 
intersection forms $E_8\oplus H$ and $E_8 \oplus 2 H$. This embedding is 
described in %Looijenga 
\cite{Lo} and \cite{FS1}. %Fintushel-Stern [FS1]. 
One can keep decomposing the obtained two manifolds as in 
\cite{FS2} %Fintushel-Stern [FS2] 
and finally split it into five pieces. Among additional embedded
homology spheres, there is 
the only one known homology 3-sphere with finite fundamental group,
the Poincar\'e homology sphere $\Sa=\Sigma(2,3,5)$. One can introduce
a spherical geometry  on $\Sa$ by representing $\pi_1(\Sa)$ as a finite
subgroup $\Ga(2,3,5)$ of the orthogonal group $O(4)$ acting free on 
$S^3=\p B^2_{\bc}$. Then $\Sa(2,3,5)=S^3/\Ga(2,3,5)$ can obtained by 
identifying the opposite sides of the spherical dodecahedron 
whose dihedral angles are $2\pi/3$, see \cite{KAG}.
\hfil\hfil\hfil \hfil\hfil\hfil\hfil\hfil\hfil\hfil\hfil\hfil
\hfil\hfil\hfil \qed
\enddemo

However we note that it is unknown whether the obtained blocks may support 
some homogeneous 4-geometries classified by Filipkiewicz \cite{F} and 
(from the point of view of Ka\"hler structures) C.T.C.Wall \cite{Wa}. %Although the boundary
This raises a question whether homogeneous geometries or splitting along homology spheres
(important from the topological point of view) are relevant for a geometrization 
of smooth 4-dimensional manifolds. For example,
neither of $Y_i$ blocks in Example 5.1 
(with the intersection form $H$) can support a complex hyperbolic structure
(which is a natural geometric candidate since $\Sa$ has a spherical CR-structure)  
because each of them has two compact boundary components. 

In fact,
in a sharp contrast to the real hyperbolic case, 
for a compact manifold $M(G)$ (that is for a geometrically finite
group $G\subset PU(n,1)$ without cusps), an application of Kohn-Rossi analytic %\cite{KR} 
extension theorem shows that the boundary of $M(G)$ is connected, and the
limit set $\La(G)$ is in some sense small (see \cite { EMM} and, 
for quaternionic and Caley hyperbolic manifolds \cite { C, CI}).
Moreover, according to a recent result of D.Burns 
(see also Theorem 4.4 in \cite{NR1}), the same claim about 
connectedness of the boundary $\p M(G)$ 
still holds if only a boundary component is compact. (In dimension $n\geq 3$,
D.Burns theorem based on \cite{BuM} uses the last compactness condition to prove
geometrical finiteness of the whole manifold $M(G)$, see also \cite{NR2}.) 

However, if no component of $\p M(G)$ is compact and we have no
finiteness condition on the holonomy group of the complex hyperbolic manifold
$M(G)$, the situation is completely different due to our construction
\cite{AX1}:

\proclaim{Theorem 5.2} In any dimension $n\geq 2$ and for any integers 
$k, k_0$, $k\geq k_0\geq 0$,
there exists a complex hyperbolic $n$-manifold $M=\ch n/G$,
$G\subset PU(n,1)$, whose boundary at infinity splits up into
$k$ connected %$(n-1)$-
manifolds, $\p_\infty M=N_1\cup\cdots\cup N_k$.
Moreover, for each boundary component $N_j$, $j\leq k_0$, its inclusion
into the manifold $M(G)$, $i_j\col N_j\subset M(G)$,
induces a homotopy equivalence of $N_j$ to $M(G)$.
\endproclaim

For a torsion free discrete group $G\subset PU(n,1)$,
a connected component $\Om_0$ of the discontinuity set $\Om(G)\subset \p\ch n$ with
the stabilizer $G_0\subset G$ is contractible and $G$-invariant if and only if 
the inclusion $N_0=\Om_0/G_0\subset M(G)$ induces a homotopy equivalence of $N_0$ to
$M(G)$ \cite{A1, AX1}. It allows us to reformulate Theorem 5.2 as

\proclaim{Theorem 5.3}  In any complex dimension $n\geq 2$ and for any natural 
numbers $k$ and $k_0$, $k\geq k_0\geq 0$, there exists a discrete
group $G=G(n,k,k_0)\subset PU(n,1)$
whose discontinuity set $\Om(G)\subset \p \ch n$ splits up into $k$
$G$-invariant components, $\Om(G)=\Om_1\cup\cdots\cup\Om_k$,
and the first $k_0$ components are contractible.
\endproclaim
 
\demo{Sketch of Proof} To prove this claim (see \cite{AX1} for details), it is crucial to construct a discrete group 
$G\subset PU(n,1)$ whose discontinuity set consists of two $G$-invariant topological
balls. To do that, we construct an infinite family
$\Sa$ of disjoint closed Heisenberg balls $B_i=B(a_i,r_i)\subset \p \ch n$ such
that the complement of their closure,
$\p \ch n \setminus \overline {\bigcup_i B(a_i,r_i)}= P_1\cup P_2$,
consists of two topological balls, $P_1$ and $P_2$. In our construction of 
such a family $\Sa$ of $\sch$-balls $B_j$,
we essentially relie on the contact structure of the Heisenberg group 
$\sch_n$. Namely, $\Sa$ is the disjoint union of finite sets $\Sa_i$
of closed $\sch$-balls whose boundary $\sch$-spheres have ``real hyperspheres"
serving as the boundaries of $(2n-2)$-dimensional cobordisms $N_i$. In the limit,
these cobordisms converge to the set of limit vertices of the polyhedra 
$P_1$ and $P_2$ which are bounded by the $\sch$-spheres $S_j=\p B_j$, 
$B_j\in\Sa$. Then the desired group $G=G(n,2,2)\subset PU(n,1)$ 
is generated by involutions $\sci_j$ which preserve those real 
$(2n-3)$-spheres lying in $S_j\subset\p P_1\cup \p P_2$, see Fig.1.

\midinsert
%\vspace{4truein}
$$\hss\psfig{file=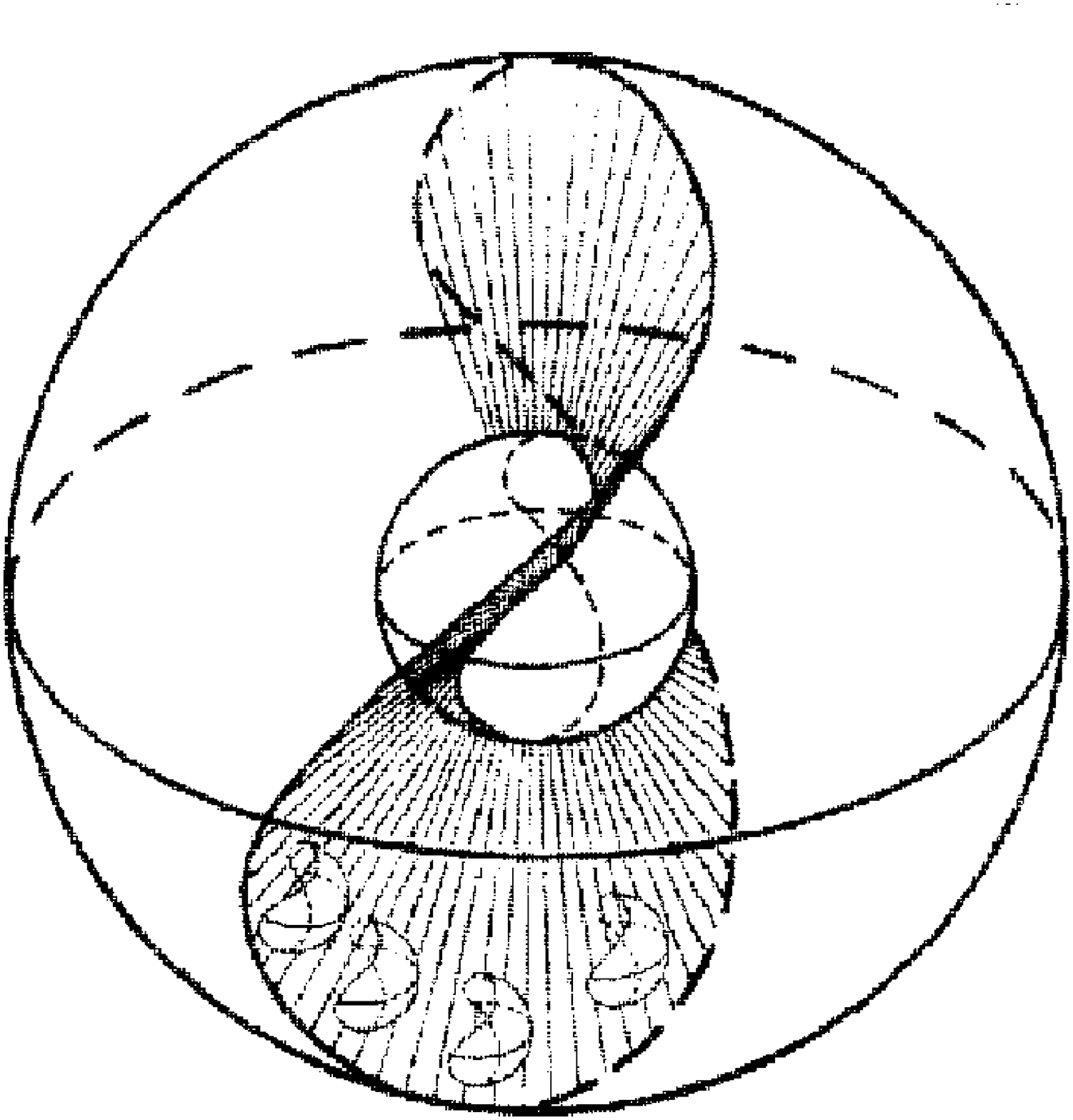,width=3truein}\hss$$
\captionwidth{32pc}
\botcaption{Figure 1} Cobordism $N_0$ in $\sch$ with two boundary real circles
\endcaption
\endinsert
\vskip5pt

We notice that, due to our construction, the intersection of each $\sch$-sphere
$S_j$ and each of the polyhedra $P_1$ and $P_2$ in the complement to the
balls $B_j\in \Sa$ is a topological $(2n-2)$-ball which splits 
into two sides,
$A_j$ and $A'_j$, and $\sci_i(A_i)=A'_i$. This allows us to define
our desired discrete group $G=G_{\Sa}\subset PU(n,1)$ as the discrete
free product,
$G_{\Sa}=\ast_j\,\, \Ga_j=\ast_i\,\, \langle \sci_j\rangle$,
of infinitely many cyclic groups $\Ga_j$
generated by involutions $\sci_j$ with respect to the $\sch$-spheres
$S_j=\p B_j$. So $P_1\cup P_2$ is a fundamental polyhedron for
the action of $G$ in $\p \ch n$, and sides of each of its connected 
components, $P_1$ or
$P_2$, are topological balls pairwise equivalent with respect to the corresponding
generators $\sci_j\in G$. Applying standard arguments (see \cite {A1}, Lemmas
3.7, 3.8), we see that the discontinuity set $\Om(G)\subset \sch_n$
consists of two $G$-invariant topological balls $\Om_1$ and $\Om_2$,
$\Om_k=\ins\left (\bigcup_{g\in G} g(\ov{P_k})\right )$,  $k=1,2$.
The fact that $\Om_k$ is a topological
ball follows from the observation that this domain is the union of a
monotone sequence,
$$
V_0=\ins (\ov{P_k})\subset V_1=\ins\left(\ov{P_k}\cup \sci_0(\ov{P_k})\right)
\subset V_2\subset\dots\,,
$$
of open topological balls, see \cite{Br}. Note that here we use the property
of our construction that $V_i$ is always a topological ball.

In the general case of $k\geq k_0\geq 0$, $k\geq 3$, we can apply the 
above infinite free products and our cobordism construction of infinite
families of $\sch$-balls with preassigned properties in order to
(sufficiently closely) ''approximate" a given hypersurface in $\sch_n$ 
by the limit sets of constructed discrete groups. For such hypersurfaces, we use the
so called ''tree-like surfaces" which are boundaries of regular neighborhoods
of trees in $\sch_n$. This allows us to
generalize A.Tetenov's \cite{T1, KAG} construction of
discrete groups $G$ on the $m$-dimensional sphere $S^m\,,\,m\geq 3$, whose
discontinuity sets split into any given number $k$ of
$G$-invariant contractible connected components.
 
\line{\hfil \hfil \hfil \qed}
\enddemo

Although, in the general case of complex hyperbolic manifolds $M$ with finitely
generated $\pi_1(M)\cong G$,
the problem on the number of boundary components of $M(G)$ is still unclear,
we show below that the situation described in Theorem 5.3 is
impossible if $M$ is geometrically finite. We refer the reader to 
\cite {AX1} for more precise formulation and proof
of this cobordism theorem: 

\proclaim {Theorem 5.4} Let $G\subset PU(n,1)$ be a geometrically
finite non-elementary torsion free discrete group whose Kleinian
manifold $M(G)$ has non-compact boundary $\p M=\Om(G)/G$ with
a component $N_0\subset \p M$ homotopy equivalent
to $M(G)$. Then there exists a
compact homology cobordism $M_c\subset M(G)$ such that $M(G)$ 
can be reconstructed from $M_c$ by gluing up a finite number of
open collars $M_i\times [0,\infty)$ where each $M_i$ is finitely covered by
the product $E_k\times B^{2n-k-1}$ of a closed (2n-1-k)-ball and a
closed $k$-manifold $E_k$ which is either flat or a 
nil-manifold (with 2-step nilpotent fundamental group).
\endproclaim

In connection to this cobordism theorem, it is worth to mention another 
interesting fact due to Livingston--Myers \cite{My} construction. Namely,
any $\bz$-homology 3-sphere is homology cobordant to a real hyperbolic
one. However, it is still unknown whether one can introduce a geometric structure
on such a homology cobordism, or a CR-structure
on a given real hyperbolic 3-manifold (in particular, a homology sphere)
or on a $\bz$-homology 3-sphere of plumbing  type. We refer to \cite{S, Mat} %A1-Sa4} and \cite{Mat} 
for recent advances on homology cobordisms, 
in particular, for results on Floer homology of homology 3-spheres and 
a new Saveliev's (presumably, homology cobordism) invariant based on 
Floer homology.

\head 6. Homeomorphisms induced by group isomorphisms
\endhead

As another application of the developed methods, we 
study the following well known problem of
geometric realizations of group isomorphisms:
% of discrete groups $G, H\subset PU(n,1)$:

\proclaim {Problem 6.1} Given a type preserving isomorphism $\varphi\col G\ra H$
of discrete groups $G, H\subset PU(n,1)$, find
subsets $X_G, X_H \subset \ov{\ch n}$ invariant for the action of groups 
$G$ and $H$, respectively,
and an equivariant homeomorphism $f_{\varphi} \col X_G\ra X_H$ which
induces the isomorphism $\varphi $. Determine metric properties of
$f_{\varphi}$, in particular, whether it is either quasisymmetric or
quasiconformal. 
\endproclaim

Such type problems were studied by several authors. In the case of
lattices $G$ and $H$ in rank 1 symmetric spaces $X$, G.Mostow \cite{Mo1}
proved in his celebrated rigidity theorem 
that such isomorphisms $\varphi\col G\ra H$ can be extended to inner isomorphisms
of $X$, provided that there is no analytic homomorphism of $X$ onto $PSL(2, \br)$. 
For that proof, it was essential to prove that $\varphi$ can be induced by
a quasiconformal homeomorphism of the sphere at infinity $\p X$  
which is the one point compactification of a (nilpotent) Carnot group $N$ 
(for quasiconformal mappings 
in Heisenberg and Carnot groups, see \cite {KR, P}).
%A.Koranyi and M.Reimann \cite {KR} and P.Pansu \cite{P}). 

If geometrically finite groups $G,H\subset PU(n,1)$ have parabolic elements and
are neither lattices nor trivial, the only 
results on geometric realization of their isomorphisms are known 
in the real hyperbolic space \cite {Tu}. Generally, those methods cannot be used in 
the complex hyperbolic space due to lack of control over convex hulls 
(where the convex hull of three points may be 4-dimensional), especially nearby
cusps. Another (dynamical) approach due to C.Yue \cite{Yu2, Cor.B} (and 
the Anosov-Smale stability theorem for hyperbolic flows) can be used 
only for convex cocompact groups $G$ and $H$, see \cite{Yu3}.
As a first step in solving the general Problem 6.1, we have 
the following isomorphism theorem \cite {A7}:

 \proclaim{Theorem 6.2} Let $\phi : G\ra H$ be a type preserving isomorphism
 of two non-ele\-men\-ta\-ry geometrically finite groups 
$G,H\subset PU(n,1)$. Then
 there exists a unique equivariant ho\-meo\-mor\-phism $f_{\phi}\col \La(G)\ra \La(H)$ of
 their limit sets that induces the isomorphism $\phi$. Moreover, 
if $\La(G)=\p \ch n$, the homeomorphism $f_{\phi}$ is the restriction  of a
hyperbolic isometry $h\in PU(n,1)$.
 \endproclaim

\demo{Proof} To prove this claim, we consider the Cayley graph $K(G,\sa)$ of a group $G$
with a given finite set $\sa$ of generators. This is a 1-complex whose vertices 
are elements of $G$, and such that two vertices $a,b\in G$ are joined by an edge
if and only if $a=bg^{\pm 1}$ for some generator $g\in \sa$. Let $|*|$ be the word norm
on $K(G,\sa)$, that is, $|g|$ equals the minimal length of words in the alphabet
$\sa$ representing a given element $g\in G$. Choosing a function $\rho$ such that
$\rho(r)=1/r^2$ for $r>0$ and $\rho(0)=1$, one can define the length
of an edge $[a,b]\subset K(G,\sa)$ as 
$d_\rho(a,b)=\min \{\rho(|a|), \rho(|b|)\}$. Considering paths of minimal length
in the sense of the function $d_\rho(a,b)$, one can extend it to a metric on 
the Cayley graph $K(G,\sa)$. So taking the Cauchy completion
$\ov {K(G,\sa)}$ of that metric space, we have the definition of the group completion 
$\ov G$ as the compact metric space 
$\ov {K(G,\sa)}\bs K(G,\sa)$, see \cite{Fl}. Up to a Lipschitz equivalence,
this definition does not depend on $\sa$. It is also clear 
that, for a cyclic group $\bz$, its completion
$\ov{\bz}$ consists of two points. Nevertheless,
for a nilpotent group $G$ with one end, its completion 
$\ov G$ is a one-point set \cite{Fl}.

Now we can define a proper equivariant embedding 
$F: K(G,\sa)\hookrightarrow \ch n$ of the Cayley graph of a given 
geometrically finite group $G\subset PU(n,1)$. To do that we may 
assume that the stabilizer of a point, say $0\in\ch n$, is trivial.
Then we set $F(g)=g(0)$ for any vertex $g\in K(G,\sa)$,
and $F$ maps any edge $[a,b]\subset K(G,\sa)$ to the geodesic segment
$[a(0),b(0)]\subset \ch n$. 

\proclaim{Proposition 6.3}  For a geometrically finite
discrete group $G\subset PU(n,1)$, there are  
constants $K,K'>0$ such that the following bounds hold for all
elements $g\in G $ with $|g|\geq K'$:
$$
 \ln(2|g|-K)^2 -\ln K^2 \leq d(0,g(0))\leq K|g|\,.\tag6.4
$$
\endproclaim

The proof of this claim is based on a comparison of the Bergman metric
$d(*,*)$ and the path metric $d_0(*,*)$ on the following subset 
$bh_0\subset \ch n$. Let $C(\La(G))\subset \ch n$ be the convex hull 
of the limit set $\La(G)\subset \p\ch n$, that is the minimal convex
subset in $\ch n$ whose closure in $\ov{\ch n}$ contains $\La(G)$.
Clearly, it is $G$-invariant, and its quotient $C(\La(G))/G$ is the minimal 
convex retract of $\ch n/G$. Since $G$ is geometrically finite, the complement
in $M(G)$ to neighbourhoods of (finitely many) cusp ends is compact
and correspond to a compact subset in the minimal convex retract, 
which can be taken as $\bh_0/G$. In other words, $\bh_0\subset C(\La(G))$
is the complement in the convex hull to a $G$-invariant family of disjoint horoballs
each of which is strictly invariant with respect to its (parabolic)
stabilizer in $G$, see \cite{AX1, Bow}, cf. also \cite{A1, Th. 6.33}.
Now, having co-compact action of the group $G$ on the domain 
$\bh_0$ whose boundary includes some horospheres, we can reduce our 
comparison  of distances $d=d(x,x')$ and $d_0=d_0(x,x')$ to their 
comparison on a horosphere. So we can take points
$x=(0,0,u)$ and $x'=(\xi,v,u)$ on a ``horizontal" horosphere 
$S_u=\bc^{n-1}\times\br\times\{u\}\subset \ch n$. Then the distances
$d$ and $d_0$ are as follows \cite{Pr2}:
$$
\cosh^2 \frac{d}{2} =\frac{1}{4u^2}\left(|\xi|^4+4u|\xi|^2+4u^2+v^2\right)\,,
 \quad d_0^2=\frac{|\xi|^2}{u} + \frac{v^2}{4u^2}\,.\tag6.5
$$

This comparison and the basic fact due to Cannon \cite{Can}
that, for a co-compact action of a group $G$ in a metric space $X$,
its Cayley graph can be quasi-isometrically embedded into $X$, finish
our proof of (6.4).

Now we apply Proposition 6.3 to define a $G$-equivariant extension of 
the map $F$ from the Cayley graph $K(G,\sa)$ to the group 
completion $\ov G$. Since the group completion of any parabolic 
subgroup $G_p\subset G$ is either a point or a two-point set
(depending on whether $G_p$ is a finite extension of cyclic or a 
nilpotent group with one end), we get 

\proclaim{Theorem 6.6} For a geometrically finite discrete group
$G\subset PU(n,1)$, there is a continuous $G$-equivariant map
$\Phi_G : \ov G\ra \La(G)$. Moreover, the map $\Phi_G$
is bijective everywhere but the set of parabolic fixed points 
$p\in\La(G)$ whose stabilizers $G_p\subset G$ have rank one. 
On this set, the map $\Phi_G$ is two-to-one.
\endproclaim

Now we can finish our proof of Theorem 6.2 by looking at the following
diagram of maps:
$$
\CD
\La(G) @<{\Phi_G}<< \ov G @>{\ov\phi}>> \ov H @>{\Phi_H}>> \La(H)\,,
\endCD
$$
where the homeomorphism $\ov\phi$ is induced by the isomorphism
$\phi$,
 and the continuous maps $\Phi_G$ and $\Phi_H$
are defined by Theorem 6.6. Namely, one can define a map
$f_\phi=\Phi_H\ov\phi\Phi^{-1}_G$. Here the map $\Phi^{-1}_G$
is the right inverse to $\Phi_G$, which exists due to Theorem 6.6.
Furthermore, the map $\Phi^{-1}_G$ is bijective 
everywhere but the set of parabolic fixed points 
$p\in\La(G)$ whose stabilizers $G_p\subset G$ have rank one, where it is 
2-to-1. Hence the composition map
$f_\phi$ is bijective and $G$-equivariant. Its uniqueness follows from 
its continuity and the fact that the image of the attractive 
fixed point of an loxodromic element $g\in G$ must be the attractive
fixed point of the loxodromic element $\phi(g)\in H$ (such loxodromic fixed points
are dense in the limit set, see \cite{A1}).

The last claim of the Theorem 6.2 directly follows from the Mostow 
rigidity theorem \cite{Mo1} because a geometrically finite group
$G\subset PU(n,1)$ with $\La(G)=\p \ch n$  is co-finite:
 $\text{\rm vol}\, (\ch n/G)<\infty$.

\line{\hfil \hfil \qed}
\enddemo

\remark{Remark 6.7}  Our proof of Theorem 6.2 can be easily
extended to the general situation, that is, to construct 
equivariant homeomorphisms $f_\phi\col \La(G)\to \La(H)$ conjugating the 
actions (on the limit sets) of isomorphic
geometrically finite groups $G, H\subset \is X$ in a (symmetric) space
$X$ with pinched negative curvature $K$, $-b^2\leq K\leq -a^2<0$.
Actually,
bounds similar to (6.4) in Prop. 6.3 (crucial for our argument) can be obtained from 
a result due to Heintze and Im Hof \cite{HI, Th.4.6} which compares the 
geometry of horospheres $S_u\subset X$ with that in the spaces of constant curvature 
$-a^2$ and $-b^2$, respectively. It gives, that for all $x,y\in S_u$ and
their distances $d=d(x,y)$ and  $d_u=d_u(x,y)$ in the space $X$
and in the horosphere $S_u$, respectively, one has that
$\frac{2}{a}\sinh (a\cdot d/2)\leq d_u \leq \frac{2}{b}\sinh (b\cdot d/2)$.
%\frac{2d}{b}\,.$$
\endremark  
\vskip4pt
 
Upon existence of such homeomorphisms $f_{\varphi}$ inducing given isomorphisms
$\varphi$ of discrete subgroups of $PU(n,1)$, the Problem 6.1
can be reduced to the questions whether $f_{\varphi}$ 
is quasisymmetric with respect to the Carnot-Carath\'eodory (or Cygan) metric, 
and whether there exists its $G$-equivariant extension to a bigger set
(to the sphere at infinity $\p X$ or even to the whole space $\ov{\ch n}$) 
inducing the isomorphism $\varphi$. For convex cocompact groups obtained by 
nearby representations, this may be seen as a generalization of
D.Sullivan stability theorem \cite{Su2}, see also \cite{A9}.

However, in a deep contrast to the real hyperbolic case,
here we have an interesting effect related to possible noncompactness of 
the boundary $\p M(G)=\Om(G)/G$. Namely, even for the simplest case of parabolic cyclic groups
$G\cong H \subset PU(n,1)$, the homeomorphic CR-manifolds 
$\p M(G)=\sch^n/G$ and $\p M(H)=\sch^n/H$ may be not quasiconformally
equivalent, see \cite {Min}. In fact, among such Cauchy-Riemannian 
3-manifolds (homeomorphic to $\br^2\times S^1$), there are exactly two quasiconformal 
equivalence classes whose representatives
have the holonomy groups generated correspondingly by a vertical $\sch$-translation
by $(0, 1)\in \bc\times \br$ and a horizontal $\sch$-translation
by $(1, 0)\in \bc\times \br$. 

Theorem 7.1 presents a more sophisticated topological
deformation $\{f_{\a}\}$, $f_{\a}\col \ov{\ch 2}\to\ov{\ch 2}$,
of a "complex-Fuchsian" co-finite group $G\subset PU(1,1)\subset PU(2,1)$ 
to quasi-Fuchsian discrete groups $G_{\a}=f_{\a}Gf_{\a}^{-1}\subset PU(2,1)$. 
It deforms pure parabolic subgroups in $G$ to subgroups in  $G_{\a}$ generated by Heisenberg 
``screw translations". 
As we point out, any such $G$-equivariant conjugations of the groups 
$G$ and $G_{\a}$ cannot be contactomorphisms because they must map some 
poli of Dirichlet bisectors to non-poli ones in the image-bisectors; 
moreover, they cannot be
quasiconformal, either. This shows the impossibility of the mentioned extension of 
Sullivan's stability theorem to the case of groups with rank one cusps.
 
Also we note that, besides the metrical (quasisymmetric) part of the
geometrization Problem 6.1, there are some topological obstructions for extensions 
of equivariant homeomorphisms $f_{\varphi}$, $f_{\varphi}\col \La(G)\ra \La(H)$. 
It follows from the next example.

\proclaim{Example 6.7} Let $G\subset PU(1,1)\subset PU(2,1)$ and 
$H\subset PO(2,1)\subset PU(2,1)$ be two geometrically finite
(loxodromic) groups isomorphic to the fundamental group $\pi_1(S_g)$ of 
a compact oriented surface $S_g$ of genus $g>1$. Then the equivariant
homeomorphism $f_{\varphi}\col \La(G)\ra \La(H)$ cannot be homeomorphically 
extended to the whole sphere $\p \ch 2 \approx S^3$.
 \endproclaim

\demo{Proof} The obstruction in this example is topological and is due to the fact that 
the quotient manifolds $M_1=\ch 2/G$ and $M_2=\ch 2/H$ are not homeomorphic. Namely, these
complex surfaces are disk bundles over the Riemann surface $S_g$ and have different Toledo
invariants: $\tau(\ch 2/G)=2g-2$ and $\tau(\ch 2/H)=0$, see \cite {To}.

The complex structures of the complex surfaces $M_1$ and $M_2$ are quite different, too.
The first manifold $M_1$ has a natural
embedding of the Riemann surface $S_g$ as a holomorphic totally geodesic closed
submanifold, and hence $M_1$ cannot be a Stein manifolds.  
The second manifolds 
$M_2$ is a Stein manifold due to a result by  Burns--Shnider \cite{BS}. 
Moreover due to Goldman \cite{Go1}, since the surface $S_p\subset M_1$ is closed, 
the manifold $M_1$ is locally rigid in the sense that 
every nearby representation $G\rightarrow PU(2,1)$
stabilizes a complex geodesic in $\ch 2$ and is conjugate to a
representation $G\rightarrow PU(1,1) \subset PU(2,1)$. In other words,
there are no non-trivial ``quasi-Fuchsian" deformations of $G$ and $M_1$. 
On the other hand, as we show in the next section (cf. Theorem 7.1), 
the second manifold $M_2$ has plentiful enough Teichm\"uller space of different 
``quasi-Fuchsian" complex hyperbolic structures.

\line{\hfil \hfil \qed}
\enddemo

\head 7. Deformations of complex hyperbolic and CR-structures:
flexibility versus rigidity
\endhead

Since any real hyperbolic $n$-manifold can be (totally geodesically) embedded to a complex
hyperbolic $n$-manifold $\ch n/G$, flexibility of the latter ones is
evident starting with hyperbolic structures on a Riemann surface of genus $g>1$, 
which form Teichm\"uller space, a complex analytic $(3g-3)$-manifold.
Strong rigidity starts in real dimension 3. Namely, 
due to the Mostow rigidity theorem \cite{M1}, hyperbolic structures of
finite volume and (real) dimension at least three are uniquely determined by
their topology, and one has no continuous deformations of them. Yet hyperbolic
3-manifolds have plentiful enough infinitesimal deformations and, according to 
Thurston's hyperbolic Dehn surgery theorem \cite{Th}, noncompact hyperbolic 
3-manifolds of finite volume can be approximated by compact hyperbolic 3-manifolds.

Also, despite their hyperbolic rigidity, real hyperbolic manifolds 
$M$ can be deformed as conformal manifolds,
or equivalently as higher-dimensional hyperbolic manifolds $M\times (0,1)$ 
of infinite volume. First such quasi-Fuchsian deformations were given by 
the author \cite{A2} and,
after Thurston's ``Mickey Mouse" example \cite{Th}, they were called bendings of 
$M$ along its totally geodesic hypersurfaces, see also \cite{A1, A2, A4-A6, JM, Ko, Su1}.
Furthermore, all these deformations are quasiconformally equivalent showing 
a rich supply of quasiconformal $G$-equivariant homeomorphisms in the 
real hyperbolic space $\rh n$. In particular, the limit set 
$\La(G)\subset\p\rh {n+1}$ deforms continuously from a round sphere 
$\p \rh n=S^{n-1}\subset S^n=\rh {n+1}$ into nondifferentiably embedded topological 
$(n-1)$-spheres quasiconformally equivalent to $S^{n-1}$. 

Contrasting to the above flexibility, ``non-real" hyperbolic manifolds
seem much more rigid. In particular, due to Pansu \cite{P}, quasiconformal maps in the
sphere at infinity of quaternionic/octionic hyperbolic spaces are
necessarily automorphisms, and thus there cannot be interesting
quasiconformal deformations of corresponding structures. Secondly, due 
to Corlette's rigidity theorem \cite{Co2}, such manifolds are even super-rigid
-- analogously to Margulis super-rigidity in higher rank \cite{MG1}.
Furthermore, complex hyperbolic manifolds share the above rigidity of  
quaternionic/octionic hyperbolic manifolds. Namely, due to the Goldman's
local rigidity theorem in dimension $n=2$ \cite{G1} and its extension for $n\geq 3$ \cite{GM}, 
every nearby discrete representation 
$\rho\col G\to PU(n,1)$ of a cocompact lattice $G\subset PU(n-1,1)$ stabilizes a complex
totally geodesic subspace $\ch {n-1}$ in $\ch n$. Thus the limit set 
$\La(\rho G)\subset \p\ch n$ is always a round sphere $S^{2n-3}$. 
In higher dimensions $n\geq 3$, this local rigidity 
of complex hyperbolic $n$-manifolds $M$ homotopy equivalent to their closed complex 
totally geodesic hypersurfaces is even global due to a recent Yue's rigidity
theorem \cite{Yu1}.
  
Our goal here is to show that, in contrast to rigidity of complex hyperbolic 
(non-Stein) manifolds $M$ from the above class, complex hyperbolic Stein manifolds $M$ 
are not rigid in general (we suspect that it is true for all Stein manifolds 
with ``big" ends at infinity). Such a flexibility has two aspects. 

First, we point out that the condition that the group $G\subset PU(n,1)$ preserves a complex
totally geodesic hyperspace in $\ch n$ is essential for local rigidity
of deformations only for co-compact lattices $G\subset PU(n-1,1)$. This
is due to the following our result \cite{ACG}:

\proclaim{Theorem 7.1} Let $G\subset PU(1,1)$ be a {\it co-finite} free lattice whose action
in $\ch 2$ is generated by four real involutions (with fixed mutually tangent 
real circles at infinity).
Then there is a continuous family $\{f_\a \}$, $-\e<\a<\e$, of $G$-equivariant homeomorphisms 
in $\ov{\ch 2}$ which induce non-trivial quasi-Fuchsian (discrete
faithful)
representations $f^*_\a \col G\to PU(2,1)$. Moreover, for each $\a\neq 0$, 
any $G$-equivariant homeomorphism of $\ov{\ch 2}$ that induces the 
representation $f^*_\a$ cannot be quasiconformal.
\endproclaim

This and an Yue's \cite{Yu2} result on Hausdorff dimension show that there are deformations
of a co-finite Fuchsian group $G\subset PU(1,1)$ into quasi-Fuchsian groups 
$G_\a=f_\a G f^{-1}_\a\subset PU(2,1)$ with Hausdorff dimension of the limit set $\La(G_\a)$
strictly bigger than one.

Secondly, we point out that the noncompactness condition in the above nonrigidity
is not essential, either. Namely, complex hyperbolic Stein manifolds $M$ 
homotopy equivalent to their closed totally 
{\it real} geodesic surfaces are not rigid, too. Namely, we give 
a canonical construction of continuous non-trivial quasi-Fuchsian
deformations of manifolds $M$, $\dim_{\bc}M=2$,
fibered over closed Riemann surfaces, which are the first such
deformations of fibrations over compact base (for a non-compact base
corresponding to an ideal triangle group
$G\subset PO(2,1)$, see \cite{GP2}).

Our construction is inspired by the approach the author used for bending deformations 
of real hyperbolic (conformal) manifolds along totally geodesic
hypersurfaces (\cite{A2, A4})
and by an example of M.Carneiro--N.Gusevskii \cite {Gu} constructing a
non-trivial discrete representation of a surface group into $PU(2,1)$.
In the case of complex hyperbolic (and Cauchy-Riemannian) structures,
the constructed ``bendings" work however in a different way than in the real case.
Namely our complex bending deformations involve
simultaneous bending of the base of the fibration of the complex surface $M$ as well 
as bendings of each of its totally geodesic fibers (see Remark 7.9). 
Such bending deformations of complex surfaces are associated
to their real simple closed geodesics (of real codimension 3), but have
nothing common with the so called cone deformations of real hyperbolic 3-manifolds along
closed geodesics, see \cite{A6, A9}. 

Furthermore, there are well known complications 
in constructing equivariant homeomorphisms in the complex hyperbolic space and 
in Cauchy-Riemannian geometry, which are due to necessary invariantness of
the K\"ahler and contact structures (correspondingly 
in $\ch n$ and at its infinity, $\ov{\sch_n}$). Despite that, the constructed
complex bending deformations are induced by equivariant 
homeomorphisms of $\ov{\ch n}$, which are in addition quasiconformal:
 
\proclaim{Theorem 7.2} Let $G\subset PO(2,1)\subset PU(2,1)$ be a given
(non-elementary) discrete group. Then, for any simple closed geodesic 
$\a$ in the Riemann 2-surface $S=H^2_{\br}/G$ and a sufficiently small $\eta_0>0$, 
there is a holomorphic family of $G$-equivariant quasiconformal homeomorphisms 
$F_{\eta}: \ov{\ch 2} \ra \ov{\ch 2}$, $-\eta_0<\eta<\eta_0$, which
defines the bending (quasi-Fuchsian) deformation 
$\Cal B_\a\col (-\eta_0,\,\eta_0)\ra \Cal R_0(G)$ of the group $G$ along the geodesic 
$\a$, $\Cal B_\a(\eta)=F^*_{\eta}$. 
\endproclaim

We notice that deformations of a complex hyperbolic manifold $M$
may depend on many parameters described by the Teichm\"uller space 
$\Cal T(M)$ of isotopy classes of complex hyperbolic structures on $M$. One can reduce
the study of this space $\Cal T(M)$ to studying the variety
$\Cal T(G)$ of conjugacy classes of discrete faithful representations 
$\rho\col G\ra PU(n,1)$ (involving the space $\Cal D(M)$ of the developing maps, see
\cite{Go2, FG}). Here $\Cal T(G)=\Cal R_0(G)/PU(n,1)$, and the variety 
$\Cal R_0(G)\subset \Hom (G, PU(n,1))$ consists of discrete faithful 
representations $\rho$ of the group $G$ with infinite co-volume, $\Vol (\ch n/G)=\infty$.
In particular, our complex bending deformations 
depend on many independent
parameters as it can be shown by applying our construction and 
\'Elie Cartan \cite {Car} angular invariant in Cauchy-Riemannian
geometry:

\proclaim{Corollary 7.3} Let $S_p=\rh 2/G$ be a closed totally real
geodesic surface of genus $p>1$ in a given complex hyperbolic surface
$M=\ch 2/G$, $G\subset PO(2,1)\subset PU(2,1)$. Then there is an
embedding $\pi\circ\Cal B\col B^{3p-3}\hra \Cal T(M)$ of a real $(3p-3)$-ball
into the Teichm\"uller space of $M$,
defined by bending deformations along disjoint closed geodesics in $M$  and
by the projection $\pi\col \Cal D(M)\ra \Cal T(M)=\Cal D(M)/PU(2,1)$ in
the development space $\Cal D(M)$.
\endproclaim

\demo{\bf Basic Construction (Proof of Theorem 7.2)} Now we start with a totally real 
geodesic surface $S=\rh 2/G$ in the complex surface
$M=\ch 2/G$, where $G\subset PO(2,1)\subset PU(2,1)$ is a
given discrete group, and fix a simple closed geodesic $\a$ on $S$. 
We may assume that the loop $\a$ is covered by a geodesic 
$A\subset \rh 2\subset \ch 2$ whose ends at infinity are $\infty$ and the origin
of the Heisenberg group $\sch=\bc \times \br$, $\ov \sch =\p \ch 2$.
Furthermore, using quasiconformal deformations of the Riemann surface $S$ 
(in the Teichm\"uller space $\Cal T(S)$, that is, by deforming the inclusion 
$G\subset PO(2,1)$ in $PO(2,1)$ by bendings along 
the loop $\a$, see Corollary 3.3 in \cite{A10}), we can assume
that the hyperbolic length of $\a$ is sufficiently small and the radius of its 
tubular neighborhood is big enough:

\proclaim{Lemma 7.4}
Let $g_\a$ be a hyperbolic element of a non-elementary discrete group 
$G\subset PO(2,1)\subset PU(2,1)$ with translation length $\ell$ 
along its axis $A\subset \rh 2$. Then any tubular neighborhood 
$U_{\da}(A)$ of the axis $A$ of radius $\da>0$  is precisely 
invariant with respect to its stabilizer $G_0\subset G$ if
$
 \sinh (\ell/4)%\left(\frac{\ell}{4}\right)
\cdot\sinh (\da/2)%\left(\frac{\da}{2}\right)
\leq 1/2$.
%\frac{1}{2}\,.$$
Furthermore, for sufficiently small $\ell$, $\ell<4\da$, the Dirichlet 
polyhedron $D_z(G)\subset \ch 2$ of the group $G$ centered at a point $z\in A$ 
has two sides $a$ and $a'$ intersecting the axis $A$ and such that 
$g_\a(a)=a'$.
\endproclaim

Then the group $G$ and its subgroups $G_0, G_1, G_2$ in the free
amalgamated (or HNN-extension) decomposition of $G$ have Dirichlet
polyhedra $D_z(G_i)\subset \ch 2$,\linebreak 
$i=0,1,2$, centered at a point $z\in
A=(0,\infty)$, whose intersections with the hyperbolic 2-plane $\rh 2$
have the shapes indicated in Figures 2-5.

\midinsert
$$\hss\psfig{file=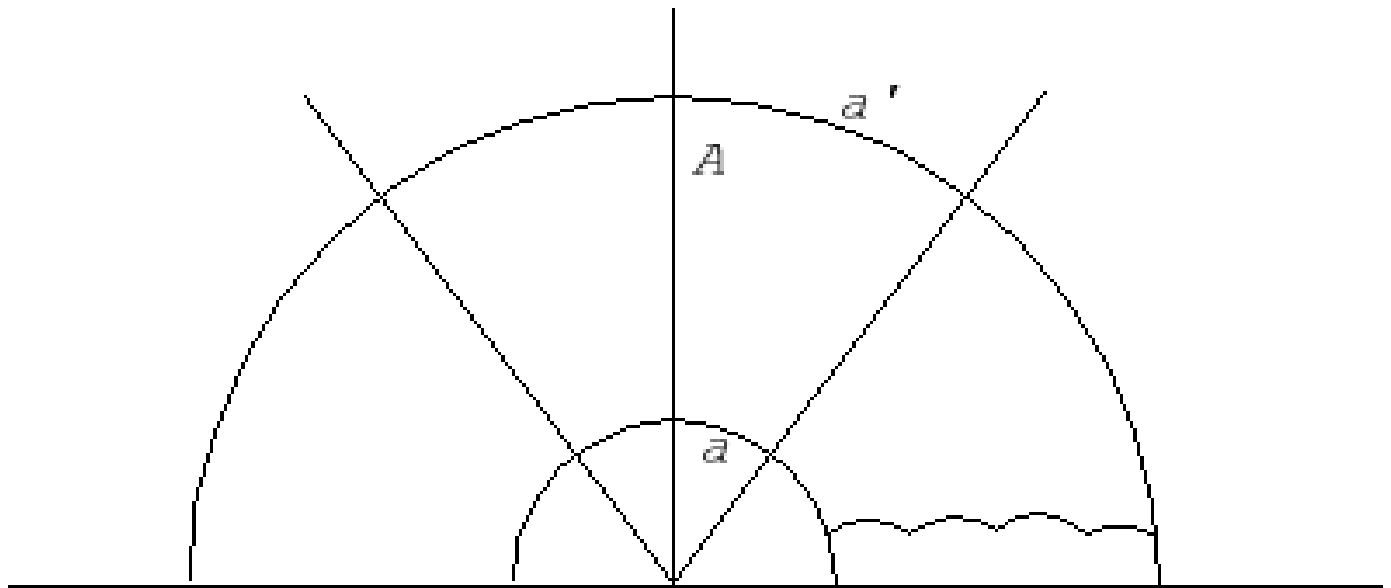,width=2.5truein}\hss \quad
\hss\psfig{file=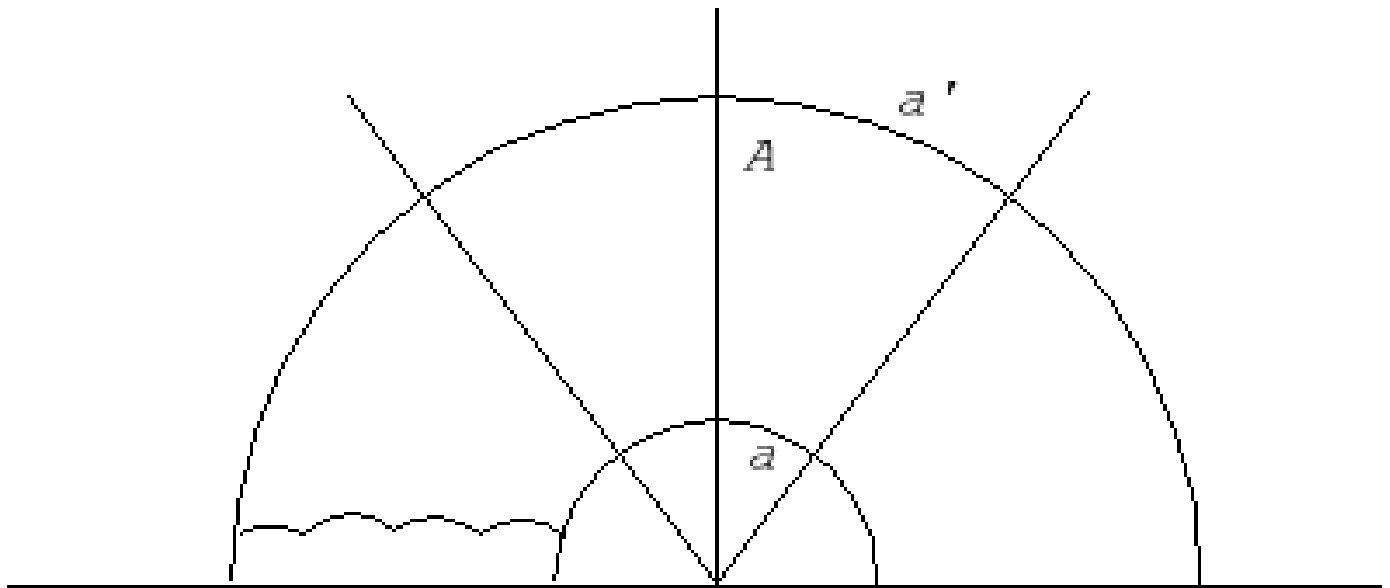,width=2.5truein}\hss$$
%\vskip1.75truein
\botcaption{Figure 2. \,\,\,
$G_1\subset G=G_1 *_{G_0}G_2$ \quad Figure 3} %\,\,\,
$G_2\subset G=G_1 *_{G_0}G_2$
\endcaption
\endinsert

%\midinsert
%$$\hss\psfig{file=bend2.eps}\hss$$
%\botcaption{Figure 3} 
%$G_2\subset G=G_1 *_{G_0}G_2$.
%\endcaption
%\endinsert

\midinsert
%\vskip1.75truein
$$\hss\psfig{file=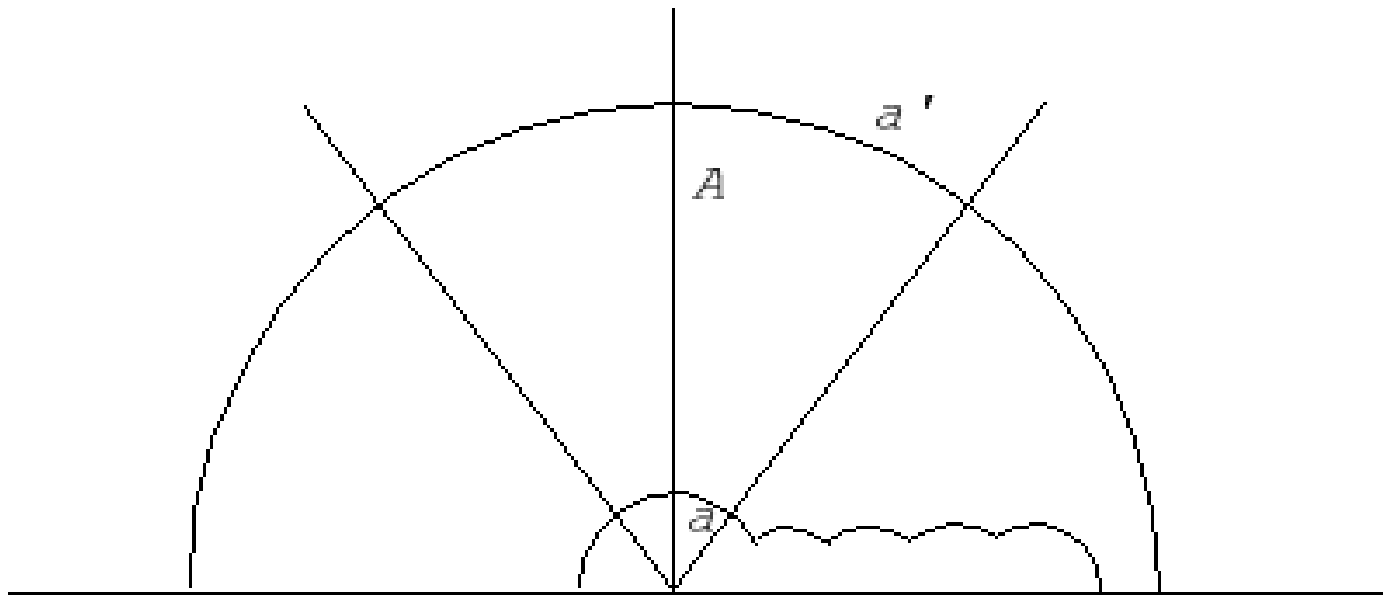,width=2.5truein}\hss\quad
\hss\psfig{file=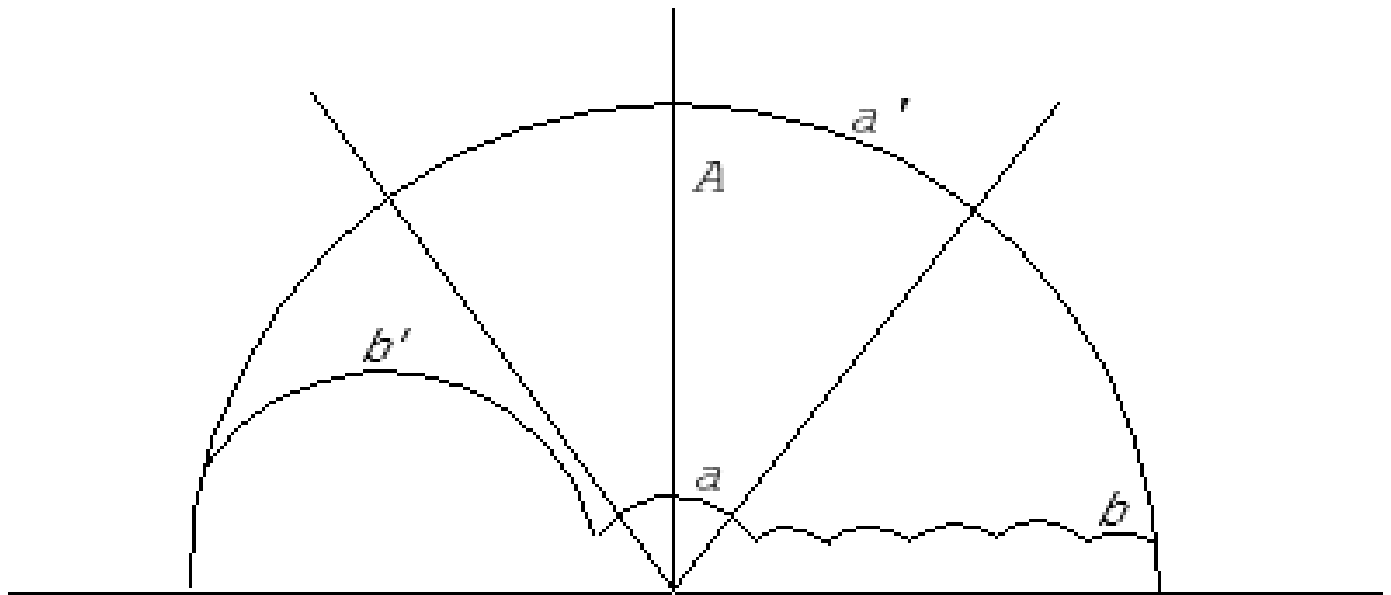,width=2.5truein}\hss$$
\botcaption {Figure 4. \,\,\, 
$G_1\subset G=G_1 *_{G_0}$ \qquad Figure 5} %\,\,\,
$G=G_1 *_{G_0}$
\endcaption
\endinsert

%\midinsert
%$$\hss\psfig{file=bend4.eps}\hss$$
%\botcaption{Figure 5} 
%$G=G_1 *_{G_0}$.
%\endcaption
%\endinsert

In particular we have that, except two bisectors $\frak S$ and $\frak S'$ that are 
equivalent under 
the hyperbolic translation $g_{\a}$ (which generates the stabilizer $G_0\subset G$ of the 
axis $A$),
all other bisectors bounding such a Dirichlet polyhedron lie in sufficiently small ``cone
neighborhoods" $C_+$ and $C_-$ of the arcs (infinite rays) $\br_+$ and $\br_-$ of the real circle
$\br\times \{0\}\subset \bc\times \br=\sch$. 

Actually, we may assume that the
Heisenberg spheres at infinity of the bisectors  $\frak S$ and $\frak S'$ 
have radii 1 and $r_0>1$, correspondingly. Then, for a
sufficiently  small $\e$, $0<\e<<r_0-1$, the cone neighborhoods 
$C_+, C_-\subset \ov {\ch 2}\bs\{\infty\}= \bc\times\br\times [0,+\infty)$
are correspondingly the cones of the $\e$-neighborhoods of the points
$(1,0,0), (-1,0,0) \in \bc\times\br\times [0,+\infty)$ with respect to the 
Cygan metric $\rho_c$ in $\ov {\ch 2}\bs\{\infty\}$, see (2.1). 

 Clearly, we may
consider the length $\ell$ of the geodesic $\a$ so small that closures of all
equidistant halfspaces in $\ov {\ch 2}\bs\{\infty\}$ bounded by those bisectors (and
whose interiors are disjoint from the Dirichlet polyhedron $D_z(G)$) do not 
intersect the co-vertical 
bisector whose infinity is $i\br\times\br\subset \bc\times \br$. It
follows from the fact \cite{Go3, Thm VII.4.0.3} that equidistant 
half-spaces  $\frak S_1$ and $\frak S_2$ in $\ch 2$ are disjoint if and
only if the half-planes $\frak S_1\cap \rh 2$ and 
$\frak S_2\cap \rh 2$ are disjoint, see Figures 2-5.

Now we are ready to define a quasiconformal bending deformation of the group $G$
along the geodesic $A$, which defines a bending deformation of the complex
surface $M=\ch 2/G$ along the given closed geodesic $\a\subset S\subset
M$.

We specify numbers $\eta$ and $\zeta$ such that
$0<\zeta<\pi/2$, $0\leq\eta<\pi-2\zeta$ and the intersection $C_+\cap (\bc\times \{0\})$
is contained in the angle $\{z\in \bc\col |\arg z|\leq\zeta\}$. Then we
define a bending homeomorphism $\phi=\phi_{\eta,\zeta}\col\bc\ra\bc$ which bends
the real axis $\br\subset\bc$ at the origin by the angle $\eta$, see Fig. 6:
$$
\phi_{\eta,\zeta}(z)=\cases
z & \text{if}\ \, |\arg z|\geq\pi-\zeta\\
z\cdot\exp(i\eta) &  \text{if}\ \, |\arg z|\leq\zeta\\
z\cdot\exp(i\eta(1-(\arg z-\zeta)/(\pi-2\zeta))) & \text{if}\ \, \zeta
<\arg z<\pi-\zeta\\
z\cdot\exp(i\eta(1+(\arg z+\zeta)/(\pi-2\zeta))) & \text{if}\ \, \zeta
-\pi<\arg z<-\zeta\,.\endcases\tag7.5
$$

\midinsert
%\vskip1.75truein %
$$\hss\psfig{file=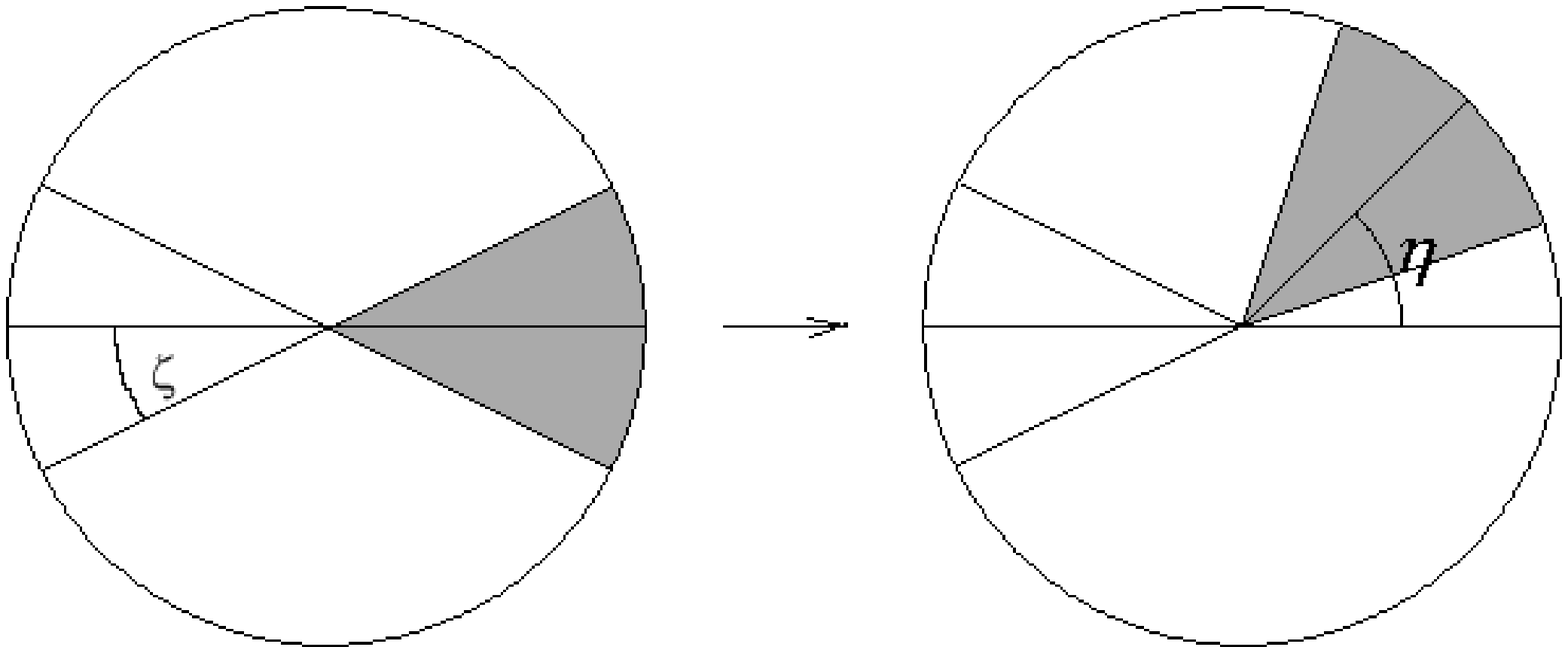,width=210pt}\hss$$
\captionwidth{7pc}
\botcaption{Figure 6} % usual picture for bending
\endcaption
\endinsert

For negative $\eta$, $2\zeta-\pi<\eta<0$, we set
$\phi_{\eta,\zeta}(z)=\ov{\phi_{-\eta,\zeta}(\ov z)}$.  Clearly,
$\phi_{\eta,\zeta}$ is quasiconformal with respect to the Cygan norm (2.1)
and is an isometry in the $\zeta$-cone
neighborhood of the real axis  $\br$ because its linear distortion is given by
$$
K(\phi_{\eta,\zeta},z)=\cases
1\, & \text{if}\ \, |\arg z|\geq\pi-\zeta\\
1\, & \text{if}\ \, |\arg z|\leq\zeta\\
(\pi-2\zeta)/(\pi-2\zeta-\eta)\, & \text{if}\ \, \zeta<\arg z<\pi-\zeta\\
(\pi-2\zeta+\eta)/(\pi-2\zeta)\, & \text{if}\ \, \zeta-\pi<\arg z<-
\zeta\,.\endcases\tag7.6
$$

Foliating the punctured Heisenberg group $\sch\bs\{0\}$ by Heisenberg
spheres $S(0,r)$ of radii $r>0$, we can extend the bending
homeomorphism $\phi_{\eta,\zeta}$ to an elementary bending homeomorphism 
$\vp=\vp_{\eta,\zeta}\col\sch\ra\sch$, $\vp(0)=0$, $\vp(\infty)=\infty$,
of the whole sphere $S^3=\ov{\sch}$ at infinity. 

Namely, for the ``dihedral angles" $W_+, W_-\subset\sch$ with the common vertical
axis $\{0\}\times \br$ and which are foliated by arcs of real circles
connecting points $(0,v)$ and $(0,-v)$ on the vertical axis and
intersecting the the $\zeta$-cone neighborhoods of infinite rays
$\br_+, \br_- \subset \bc$, correspondingly, the restrictions $\vp|_{W_-}$ 
and $\vp|_{W_+}$ of the bending homeomorphism 
$\vp=\vp_{\eta,\zeta}$ are correspondingly the identity and the unitary rotation
$U_{\eta}\in PU(2,1)$ by angle $\eta$ about the vertical axis
$\{0\}\times\br\subset \sch$, see also \cite{A10, (4.4)}.
Then it follows from 
(7.6) that $\vp_{\eta,\zeta}$ is a $G_0$-equivariant quasiconformal homeomorphism 
in $\ov{\sch}$. 

We can naturally
extend the foliation of the punctured Heisenberg
group $\sch\bs\{0\}$ by Heisenberg spheres $S(0,r)$ to a foliation of the hyperbolic 
space $\ch 2$ by bisectors $\gtS_r$ having those $S(0,r)$ as the spheres at
infinity. It is well known (see \cite{M2}) that each bisector $\gtS_r$
contains a geodesic $\ga_r$ which connects points $(0,-r^2)$ and
$(0,r^2)$ of the Heisenberg group $\sch$ at infinity, and furthermore 
$\gtS_r$ fibers over  $\ga_r$  by complex geodesics $Y$ whose circles at
infinity are complex circles foliating the sphere $S(0,r)$. 

Using those foliations of the hyperbolic 
space $\ch 2$ and bisectors $\gtS_r$, we extend the elementary bending homeomorphism 
$\vp_{\eta,\zeta}\col\ov{\sch}\ra\ov{\sch}$ at infinity to an elementary
bending homeomorphism $\vP_{\eta,\zeta}\col\ov{\ch 2}\ra\ov{\ch 2}$.
Namely, the map $\vP_{\eta,\zeta}$ preserves each of bisectors $\gtS_r$,
each complex geodesic fiber $Y$ in such bisectors, and fixes the intersection points
$y$ of those complex geodesic fibers and the complex geodesic connecting the origin and 
$\infty$ of the Heisenberg group $\sch$ at infinity. We complete our extension
$\vP_{\eta,\zeta}$ by defining its restriction to a given (invariant) complex
geodesic fiber $Y$ with the fixed point $y\in Y$. This map is obtained by  
radiating the circle homeomorphism $\vp_{\eta,\zeta}|_{\p Y}$
to the whole (Poincar\'e) hyperbolic 2-plane $Y$ along geodesic rays
$[y,\infty)\subset Y$, so that it preserves circles in $Y$ centered at $y$ and 
bends (at $y$, by the angle $\eta$) the geodesic in $Y$ connecting the central
points of the corresponding arcs of the complex circle $\p Y$, see
Fig.6.

Due to the construction, the elementary bending (quasiconformal) 
homeomorphism $\vP_{\eta,\zeta}$ commutes with elements of the cyclic loxodromic
group $G_0\subset G$. Another most important property of the 
homeomorphism $\vP_{\eta,\zeta}$ is the following. 

Let $D_z(G)$ be 
the Dirichlet fundamental polyhedron of the group $G$ centered at a given
point $z$ on the axis $A$ of the cyclic loxodromic group $G_0\subset G$,
and $\gtS^+\subset \ch 2$ be a ``half-space" disjoint from  $D_z(G)$ and 
bounded by a bisector $\gtS\subset \ch 2$ which is different from bisectors 
$\gtS_r, r>0$, and contains a side $\gts$ of the polyhedron $D_z(G)$. Then there is
an open neighborhood $U(\ov{\gtS^+})\subset \ov{\ch 2}$ such that the restriction
of the elementary bending homeomorphism $\vP_{\eta,\zeta}$ to it either
is the identity or coincides with the unitary rotation $U_{\eta}\subset PU(2,1)$
by the angle $\eta$ about the ``vertical" complex geodesic (containing 
the vertical axis $\{0\}\times\br\subset \sch$ at infinity).

The above properties of quasiconformal homeomorphism $\vP=\vP_{\eta,\zeta}$
show that the image $D_{\eta}=\vP_{\eta,\zeta}(D_z(G))$ is a polyhedron in 
$\ch 2$ bounded by bisectors. Furthermore, there is a natural identification of
its sides induced by $\vP_{\eta,\zeta}$. Namely, the pairs of sides preserved by
$\vP$ are identified by the original generators of the group $G_1\subset G$. 
For other sides $\gts_{\eta}$ of $D_{\eta}$, which are images of corresponding 
sides $\gts\subset D_z(G)$
under the unitary rotation $U_{\eta}$, we define side pairings by using the group
$G$ decomposition (see Fig. 2-5). 

Actually, if $G=G_1 *_{G_0}G_2$,
we change the original side pairings $g\in G_2$ of $D_z(G)$-sides to the 
hyperbolic isometries $U_{\eta}gU_{\eta}^{-1}\in PU(2,1)$. In the case
of HNN-extension, $G=G_1 *_{G_0}=\langle G_1, g_2\rangle$, 
we change the original side pairing $g_2\in G$ of $D_z(G)$-sides to the
hyperbolic isometry $U_{\eta}g_2\in PU(2,1)$. In other words, we define
deformed groups $G_{\eta}\subset PU(2,1)$  correspondingly as
$$
G_{\eta}=G_1 *_{G_0}U_{\eta}G_2U^{-1}_{\eta}\quad \text{or}\quad
    G_{\eta}=\langle G_1, U_{\eta}g_2\rangle=G_1 *_{G_0}\,.\tag7.7
$$
This shows that the family of representations $G\ra G_{\eta}\subset PU(2,1)$
does not depend on 
angles $\zeta$ and holomorphically depends on the angle parameter $\eta$. 
Let us also observe that, for small enough angles $\eta$, the behavior of neighboring 
polyhedra $g'(D_{\eta})$,
$g'\in G_{\eta}$ is the same as of those $g(D_z(G))$, $g\in G$, around
the Dirichlet fundamental polyhedron $D_z(G)$. This is because the new polyhedron 
$D_{\eta}\subset \ch 2$ has isometrically the same (tesselations of) 
neighborhoods of its side-intersections 
as $D_z(G)$ had. This implies that the polyhedra $g'(D_{\eta})$,
$g'\in G_{\eta}$, form a tesselation of $\ch 2$ (with non-overlapping interiors). 
Hence the deformed group $G_{\eta}\subset PU(2,1)$ is a discrete group,
and $D_{\eta}$ is its fundamental polyhedron bounded by bisectors.

Using $G$-compatibility of the restriction of the
elementary bending homeomorphism $\vP=\vP_{\eta,\zeta}$ to the closure
$\ov{D_z(G)}\subset \ov{\ch 2}$, we equivariantly extend it from the polyhedron 
$\ov{D_z(G)}$ to the whole space $\ch 2\cup\Om(G)$ accordingly to the $G$-action.
In fact, in terms of the natural isomorphism $\chi\col G\ra G_{\eta}$ which is 
identical on the subgroup $G_1\subset G$, we can write the obtained 
$G$-equivariant homeomorphism 
$F=F_{\eta}\,\col\, \ov{\ch 2}\bs\La(G)\ra\ov{\ch 2}\bs\La(G_{\eta})$ 
in the following form:
$$
\aligned F_{\eta}(x)&=\vP_{\eta}(x)\quad\text{for }\quad x\in
\ov{D_z(G)},\\
     F_{\eta}\circ g(x)&=g_{\eta}\circ F_{\eta}(x)
\quad\text{for }\quad 
x\in \ov{\ch 2}\bs\La(G),\,\, g\in G,\,\, g_{\eta}=\chi(g)\in G_{\eta}\,.
\endaligned\tag7.8
$$

Due to quasiconformality of $\vP_{\eta}$, the extended $G$-equivariant 
homeomorphism $F_{\eta}$
is quasiconformal. Furthermore, its extension by continuity to the limit
(real) circle $\La(G)$ coincides with the canonical 
equivariant homeomorphism
$f_{\chi}\col \La(G)\ra\La(G_{\eta})$ given by the isomorphism 
Theorem 6.2. Hence we have a $G$-equivariant
quasiconformal self-homeomorphism of the whole space $\ov{\ch 2}$, which we
denote as before by $F_{\eta}$. 

The family of $G$-equivariant
quasiconformal homeomorphisms $F_{\eta}$ induces representations
$F^*_{\eta}\,\col\, G\to G_{\eta}=F_{\eta}G_2F^{-1}_{\eta}$, 
$\eta\in (-\eta_0,\,\eta_0)$. In other words, we have
a curve $\Cal B\col (-\eta_0,\,\eta_0)\ra \Cal R_0(G)$
in the variety $\Cal R_0(G)$ of faithful discrete representations of $G$ 
into $PU(2,1)$, which covers a nontrivial curve
in the Teichm\"uller space $\Cal T(G)$ represented by conjugacy classes
$[\Cal B(\eta)]=[F^*_{\eta}]$. We call the constructed  deformation $\Cal B$ 
the bending deformation
of a given lattice $G\subset PO(2,1)\subset PU(2,1)$ along a bending geodesic 
$A\subset \ch 2$ with loxodromic stabilizer $G_0\subset G$. In terms of
manifolds, $\Cal B$ is the bending deformation of a given complex surface
$M=\ch 2/G$ homotopy equivalent to its totally real geodesic surface 
$S_g\subset M$, along a given simple geodesic $\a$. 

 \line{\hfil \hfil \qed}
\enddemo

\remark{Remark 7.9} It follows from the above construction of the 
bending homeomorphism $F_{\eta, \zeta}$, that the deformed complex
hyperbolic surface $M_{\eta}=\ch 2/G_{\eta}$ fibers over the pleated 
hyperbolic surface $S_{\eta}=F_{\eta}(\rh 2)/G_{\eta}$ (with the closed geodesic 
$\a$ as the singular locus). The fibers of this fibration are
``singular real planes" obtained from totally real geodesic 2-planes by
bending them  by angle $\eta$ along complete real geodesics. These
(singular) real geodesics are the intersections of the complex geodesic
connecting the axis $A$ of the cyclic group $G_0\subset G$ and the
totally real geodesic planes that represent fibers of the original fibration in
$M=\ch 2/G$.
\endremark

\demo{Proof of Corollary 7.3} Since, due to (7.7), bendings along disjoint closed geodesics
are independent, we need to show that 
our bending deformation is not trivial, and 
$[\Cal B(\eta)]\neq [\Cal B(\eta')]$ for any $\eta\neq\eta'$. 

The non-triviality of our deformation follows directly from (7.7), cf. \cite{A9}.
Namely, the restrictions $\rho_{\eta}|_{G_1}$ of bending representations
to a non-elementary subgroup $G_1\subset G$ (in general, to a ``real" subgroup
$G_r\subset G$ corresponding to a totally real geodesic piece in the homotopy
equivalent surface $S\backsimeq M$) are identical. So if the deformation $\Cal B$ 
were trivial then it would be conjugation of the group $G$ by projective
transformations that commute with the non-trivial real subgroup $G_r\subset G$ and 
pointwise fix the totally real geodesic plane $\rh 2$. This contradicts
to the fact that the limit set of any deformed group $G_{\eta}$,
$\eta\neq 0$, does not belong to the real circle containing the 
limit Cantor set $\La(G_r)$.

The injectivity of the map $\Cal B$ can be obtained by using \'Elie Cartan \cite{Car}
angular invariant  $\ba(x)$, $-\pi/2\leq\ba(x)\leq\pi/2$, 
for a triple $x=(x^0, x^1, x^2)$ of points in 
$\p\ch 2$. It is known (see \cite{Go3}) that, for two triples $x$ and $y$, 
$\ba(x)=\ba(y)$ if and only if
there exists $g\in PU(2,1)$ such that $y=g(x)$; furthermore, such a $g$ is unique
provided that $\ba(x)$ is neither zero nor  $\,\pm \pi/2$. Here 
$\ba(x)=0$ if and only if  $x^0, x^1$ and $x^2$ lie on an $\br$-circle, and
$\ba(x)=\pm \pi/2$ if and only if $x^0, x^1$ and $x^2$ lie on a
chain ($\bc$-circle).  

 Namely, let $g_2\in G\bs G_1$ be a generator of
the group $G$ in (4.5) whose fixed point $x^2\in \La(G)$ lies in $\br_+\times \{0\}
\subset \sch$, and $x^2_{\eta}\in \La(G_{\eta})$ the corresponding fixed point
of the element $\chi_{\eta}(g_2)\in G_{\eta}$ under the free-product
isomorphism $\chi_{\eta}\col G\ra G_{\eta}$.
Due to our construction, one can see that the orbit
$\ga(x^2_{\eta})$, $\ga\in G_0$, under the loxodromic (dilation) subgroup $G_0\subset 
G\cap G_{\eta}$ approximates the origin along a ray 
$(0,\infty)$ which has a non-zero
angle $\eta$ with the ray $\br_-\times \{0\}\subset \sch$. The latter ray also contains
an orbit $\ga(x^1)$, $\ga\in G_0$, of a limit point $x^1$ of $G_1$ which approximates 
the origin from the other side. Taking triples $x=(x^1, 0, x^2)$ and 
$x_{\eta}=(x^1, 0, x^2_{\eta})$ of points which lie correspondingly in
the limit sets $\La(G)$ and $\La(G_{\eta})$, we have
that $\ba(x)=0$ and $\ba(x_{\eta})\neq 0, \,\pm\pi/2$. Due to Theorem 6.2,
both limit sets are topological circles which however cannot be equivalent
under a hyperbolic isometry because of different Cartan invariants
(and hence, again, our deformation is not
trivial). 

Similarly, for two different values $\eta$ and $\eta'$, we have triples
$x_{\eta}$ and $x_{\eta'}$ with different (non-trivial) Cartan angular invariants 
$\ba(x_{\eta})\neq \ba(x_{\eta'})$. Hence 
$\La(G_{\eta})$ and $\La(G_{\eta'})$ are not $PU(2,1)$-equivalent. 

 \line{\hfil \hfil \qed}
\enddemo

One can apply the above proof to a general situation of bending
deformations of a complex hyperbolic surface $M=\ch 2/G$ whose holonomy 
group  $G\subset PU(2,1)$ has a non-elementary subgroup $G_r$ preserving
a totally real geodesic plane $\rh 2$. In other words, such a complex surfaces $M$
has an embedded totally real geodesic surface with geodesic boundary.
 In particular all complex
surfaces constructed in \cite{GKL} with a given Toledo invariant lie in
this class. So we immediately have:

\proclaim{Corollary 7.10} Let $M=\ch 2/G$ be a complex hyperbolic surface
with embedded totally real geodesic surface $S_r\subset M$
with geodesic boundary, and $\Cal B\col (-\eta, \eta)\ra \Cal D(M)$ be
the bending deformation of $M$ along a simple closed geodesic $\a\subset S_r$. 
Then the map $\pi\circ\Cal B\col (-\eta, \eta)\ra \Cal T(M)=\Cal D(M)/PU(2,1)$ is a smooth
embedding provided that the limit set $\La(G)$ of the holonomy group $G$
does not belong to the $G$-orbit of the real circle $S^1_{\br}$ and the
chain $S^1_{\bc}$, where the latter is the infinity of the complex geodesic containing 
a lift $\tilde\a\subset\ch 2$ of the closed geodesic $\a$, and the former one contains
the limit set of the holonomy group $G_r\subset G$ of the geodesic
surface $S_r$.
 
\line{\hfil\hfil\qed}
\endproclaim

As an application of the constructed bending deformations, we
answer a well known question about cusp groups on the boundary of
the Teichm\"uller space $\sct (M)$ of a Stein complex hyperbolic surface $M$
fibering over a compact Riemann surface of genus $p>1$. It is a direct corollary of the 
following result, see \cite{AG}:
  
\proclaim {Theorem 7.11} Let $G\subset PO(2,1)\subset PU(2,1)$ be a
non-elementary discrete group
$S_p$ of genus $p\geq 2$. Then, for any simple closed geodesic
$\a$ in the Riemann surface $S=H^2_{\br}/G$, there is a continuous deformation $\rho_t=f^*_t$ 
induced by $G$-equivariant quasiconformal homeomorphisms 
$f_t: \ov{\ch 2} \to \ov{\ch 2}$ whose limit representation $\rho_{\infty}$ corresponds
to a boundary cusp point of the Teichm\"uller space $\Cal T(G)$, that
is, the boundary group $\rho_{\infty}(G)$ has an accidental parabolic element 
$\rho_{\infty}(g_{\a})$ where $g_{\a}\in G$ represents the geodesic $\a\subset S$.
\endproclaim

We note that, due to our construction of such continuous quasiconformal deformations
in \cite {AG},
they are independent if the corresponding geodesics $\a_i\subset S_p$
are disjoint. It implies the existence of a boundary group in $\p \Cal T(G)$ with
``maximal" number of non-conjugate accidental parabolic subgroups:

\proclaim {Corollary 7.12} Let $G\subset PO(2,1)\subset PU(2,1)$ be a
uniform lattice isomorphic to the fundamental group of a closed surface
$S_p$ of genus $p\geq 2$. Then there is a continuous deformation 
$R\col \br^{3p-3}\to \Cal T(G)$ whose boundary group
$G_{\infty}=R(\infty)(G)$ has $(3p-3)$ non-conjugate accidental parabolic
subgroups.
\endproclaim

Finally, we mention another aspect of the intrigue Problem 4.12 on geometrical finiteness
of complex hyperbolic surfaces (see \cite {AX1, AX2}) for which it may 
perhaps be possible to 
apply our complex bending deformations:

\proclaim {Problem} Construct a geometrically infinite (finitely
generated) discrete group $G\subset PU(2,1)$ whose limit set is the whole sphere
at infinity, $\La(G)=\p \ch 2=\ov {\sch}$, and which is the limit of
convex cocompact groups $G_i\subset PU(2,1)$ from the Teichm\"uller space 
$\Cal T(\Ga)$ of a convex cocompact group $\Ga\subset PU(2,1)$. Is that
possible for a Schottky group $\Ga$?
\endproclaim
%\vfil
%\eject

\def\ref#1{[#1]}
\eightpoint
\parindent=36pt

\head  REFERENCES \endhead
\bigskip

\frenchspacing

\item{\ref{Ah}} Lars V. Ahlfors, Fundamental polyhedra and limit point sets
of Kleinian groups. - Proc. Nat. Acad. Sci. USA, {\bf 55}(1966), 251-254.

\item{\ref{A1}}  Boris Apanasov, Discrete groups in Space and
Uniformization  Problems. - Math. and Appl., {\bf 40}, Kluwer
Academic Publishers, Dordrecht, 1991.

\item{\ref{A2}}  \underbar{\phantom{apana}}, Nontriviality  of Teichm\"uller space for
Kleinian  group in  space.-  Riemann Surfaces and Related
Topics: Proc. 1978 Stony  Brook Conference (I.Kra and
B.Maskit, eds), Ann. of  Math. Studies {\bf 97},
Princeton Univ. Press, 1981, 21-31.

\item{\ref{A3}}\underbar{\phantom{xxxxx}}, Geometrically finite hyperbolic
structures on manifolds. - Ann. of Glob. Analysis and Geom., {\bf 1:3}(1983),
1-22.

\item{\ref{A4}} \underbar{\phantom{xxxxx}}, Thurston's bends and geometric
deformations of conformal structures.- Complex Analysis and
Applications'85,  Publ. Bulgarian Acad. Sci.,
Sofia, 1986, 14-28.

\item{\ref{A5}} \underbar{\phantom{xxxxx}}, Nonstandard uniformized conformal
structures on hyperbolic manifolds. - Invent. Math., {\bf 105}(1991), 137-152.

\item{\ref{A6}}  \underbar{\phantom{xxxxx}}, Deformations of conformal structures on
hyperbolic manifolds.- J. Diff. Geom. {\bf 35}
(1992), 1-20.

\item{\ref{A7}} \underbar{\phantom{xxxxx}}, Canonical homeomorphisms in
Heisenberg group induced by isomorphisms of discrete subgroups 
of $PU(n,1)$.- Preprint,1995; Russian Acad.Sci.Dokl.Math., to appear. 

\item{\ref{A8}} \underbar{\phantom{xxxxx}}, Quasiconformality and
geometrical finiteness in Carnot-Carath\'eodory and negatively curved
spaces. - Math. Sci. Res. Inst. at Berkeley, 1996-019.

\item{\ref{A9}}\underbar{\phantom{xxxxx}}, Conformal geometry of
discrete groups and manifolds.  - W. de Gruyter, Berlin- New York, to appear.

\item{\ref{A10}} \underbar{\phantom{xxxxx}},
Bending deformations of complex hyperbolic surfaces. - 
Preprint 1996-062, Math. Sci. Res. Inst., Berkeley, 1996.

\item{\ref{ACG}} Boris Apanasov, Mario Carneiro and Nikolay Gusevskii,
Some deformations of complex hyperbolic surfaces. - In preparation.
 
\item{\ref{AG}} Boris Apanasov and Nikolay Gusevskii, 
The boundary of Teichm\"uller space of complex hyperbolic surfaces. - In preparation.

\item{\ref{AX1}}  Boris Apanasov and Xiangdong Xie, Geometrically finite 
complex hyperbolic manifolds.- Preprint, 1995.

\item{\ref{AX2}} \underbar{\phantom{xxxxx}}, Manifolds of negative
curvature and nilpotent groups.- Preprint, 1995.

\item{\ref{AT}}  Boris Apanasov and Andrew Tetenov, Nontrivial cobordisms with
geometrically finite hyperbolic structures. - J. Diff. Geom., {\bf 28}(1988),
407-422.

\item{\ref{BGS}} W.~Ballmann, M.~Gromov and V.~Schroeder,  Manifolds of
Nonpositive Curvature. -\break Birkh\"auser, Boston, 1985.

\item{\ref{BM}} A.F.~Beardon and B. Maskit, Limit points of Kleinian
groups and finite sided fundamental polyhedra. - Acta Math. {\bf 132}
(1974), 1-12.

%\item{\ref{B}} A.F.~Beardon, The geometry of discrete groups.- Springer-Verlag,
%Berlin/New York, 1983.

\item{\ref{Be}} Igor Belegradek, Discrete surface groups actions with
accidental parabolics on complex hyperbolic plane.- Preprint, Univ. of
Maryland.

\item{\ref{BE}} \v Zarko Bi\v zaca and John Etnyre, Smooth structures on
collarable ends of 4-manifolds.- Preprint, 1996.

\item{\ref{Bor}} A.Borel, Compact Clifford-Klein forms of symmetric
spaces. - Topol.{\bf 2}(1962), 111-122.

\item{\ref{Bow}} Brian Bowditch,
Geometrical finiteness with variable negative curvature. -
Duke J. Math., {\bf 77} (1995), 229-274.

\item{\ref{Br}} M.~Brown, The monotone union of open $n$-cells is an open
$n$-cell. - Proc. Amer. Math. Soc. {\bf 12} (1961), 
812-814.

\item{\ref{BS}} D.~Burns and S.~Shnider, Spherical hypersurfaces in complex
manifolds. - Inv. Math., {\bf 33}(1976), 223-246.

\item{\ref{BuM}} Dan Burns and Rafe Mazzeo, On the geometry of cusps for
$SU(n,1)$. - Preprint, Univ. of Michigan, 1994.

\item{\ref{BK}} P.Buser and H.Karcher, Gromov's almost flat manifolds. -
Asterisque {\bf 81} (1981), 1-148.

\item{\ref{Can}} James W.~Cannon, The combinatorial structure of cocompact 
discrete hyperbolic groups. - Geom. Dedicata {\bf 16} (1984), 123--148.

\item{\ref{Car}} \'E.Cartan, Sur le groupe de la g\'eom\'etrie
hypersph\'erique.- Comm.Math.Helv.{\bf 4}(1932), 158%-171.

\item{\ref{CG}} S.~Chen and L.~Greenberg, Hyperbolic spaces.-
Contributions to Analysis, Academic Press, New York, 1974, 49-87.

\item{\ref{Co1}} Kevin Corlette, Hausdorff dimensions of limit 
sets. - Invent. Math. {\bf 102} (1990), 521-542.

\item{\ref{Co2}}  \underbar{\phantom{xxxxx}}, 
Archimedian superrigidity and hyperbolic geometry. - 
Ann. of Math. {\bf 135} (1992), 165-182.

\item{\ref{Cy}} J.~Cygan, Wiener's test for Brownian motion on the Heisenberg
group. - Colloquium Math. {\bf 39}(1978), 367-373.

\item{\ref{EMM}} C.Epstein, R.Melrose and G.Mendoza, Resolvent of the Laplacian
on strictly pseudoconvex domains. - Acta Math. {\bf 167}(1991), 1-106.

%\item{\ref{E}} D.B.A. Epstein, Complex hyperbolic geometry.- London
%Math. Soc. Lect. Notes {\bf 111}, Cambridge Univ. Press, 1987, 93-117.

\item{\ref{F}} R.P. Filipkiewicz, Four-dimensional geometries. -
Ph.D.Thesis, Univ. of Warwick, 1984.

\item{\ref{FG}} Elisha Falbel and Nikolay Gusevskii, Spherical CR-manifolds of
dimension 3.- Bol. Soc. Bras. Mat. {\bf 25} (1994), 31-56.

\item{\ref{FS1}} Ronald Fintushel and Ronald Stern, Homotopy $K3$ surfaces containing 
$\Sigma (2,3,7)$. - J. Diff. Geom.
{\bf 34} (1991), 255--265.

\item{\ref{FS2}} \underbar{\phantom{xxxxx}}, Using Floer's exact triangle 
to compute Donaldson invariants. -
The Floer memorial volume,
Progr. Math. {\bf 133}, Birkh\"auser, Basel, 1995, 435--444.

\item{\ref{FS3}} \underbar{\phantom{xxxxx}}, Surgery in cusp
neighborhoods and the geography of irreducible 4-manifolds.- Invent.
Math. {\bf 117} (1994), 455-523.

\item{\ref{FT}}  M.H. Freedman and  L. Taylor, $\Lambda$-splitting 
4-manifolds. -
Topology {\bf 16} (1977), 181-184.

\item{\ref{Go1}} William Goldman, Representations of fundamental groups of
surfaces.- Geometry and topology, %(J.Alexander and J.Harer, Eds), 
Lect. Notes Math. {\bf 1167}, Springer, 1985, 95-117. 

\item{\ref{Go2}} \underbar{\phantom{xxxxx}}, Geometric structures on
manifolds and varieties of representations. - Geometry of Group
Representations, %(A.Magid and W. Goldman, Eds), 
Contemp. Math. {\bf 74} (1988), 169-198.
 
\item{\ref{Go3}} \underbar{\phantom{xxxxx}}, Complex hyperbolic
geometry. - Oxford Univ. Press, to appear.

\item{\ref{GKL}} W.Goldman, M.Kapovich and B.Leeb, Complex hyperbolic
manifolds homotopy equivalent to a Riemann surface. - Preprint, 1995.

\item{\ref{GM}} William Goldman and John Millson, Local rigidity of discrete
groups acting on complex hyperbolic space.- Invent. Math. {\bf 88}(1987),
495-520.

\item{\ref{GP1}} William Goldman and John Parker, Dirichlet polyhedron for
dihedral groups acting on complex hyperbolic space. -
J. of Geom. Analysis, {\bf 2:6} (1992), 517-554. 

\item{\ref{GP2}}  \underbar{\phantom{xxxxx}}, %William Goldman and John Parker, 
Complex hyperbolic ideal
triangle groups.- J.reine angew.Math.{\bf 425}(1992),71-86.

\item{\ref{Gr}} M. Gromov, Almost flat manifolds. - J. Diff. Geom.
{\bf 13} (1978), 231-141.

\item{\ref{Gu}} Nikolay Gusevskii, Colloquium talk, Univ. of Oklahoma,
Norman, December 1995.

\item{\ref{HI}} Ernst Heintze and Hans-Christoph Im Hof, Geometry of
horospheres. - J. Diff. Geom. {\bf 12} (1977), 481--491.

\item{\ref{JM}}  Dennis Johnson and  John Millson,  Deformation  spaces
associated  to  compact  hyperbolic  manifolds.- Discrete
Groups in Geometry  and Analysis,%:  Papers in  Honor of G.D.
%Mostow on  His Sixtieth Birthday,  Ed. R. Howe, 
Birkhauser, %Boston, 
1987, 48-106.

\item{\ref{Kob}} Shoshichi Kobayashi, Hyperbolic manifolds and
holomorphic mappings.- M.Decker, 1970.

\item{\ref{KR}} Adam Koranyi and Martin Reimann, Quasiconformal mappings on the
Heisenberg group. - Invent. Math. {\bf 80} (1985), 309-338.

%\item{\ref{KR2}} \underbar{\phantom{xxxxx}}, Foundations for the
%theory of quasiconformal mappings on the Hei\-sen\-berg group. - Adv. Math.
%{\bf 111} (1995), 1-87.

\item{\ref{Ko}} Christos Kourouniotis,
Deformations of hyperbolic structures.- Math. Proc. Cambr. Phil. Soc.
 {\bf 98} (1985), 247-261.

\item{\ref{KAG}} Samuil Krushkal', Boris Apanasov and Nikolai Gusevskii,
Kleinian groups and uniformization in examples and problems. - Trans. Math.
Mono. {\bf 62}, Amer. Math. Soc., Providence, 1986, 1992.

\item{\ref{LB}} Claude LeBrun, Einstein metrics and Mostow rigidity. -
Preprint, Stony Brook, 1994.

\item{\ref{Lo}} Eduard Looijenga,
The smoothing components of a triangle singularity. -
 Proc. Symp. Pure Math.
{\bf 40} (1983 ), 173--184.

\item{\ref{Mar}} Albert Marden, The geometry of finitely generated Kleinian
groups. - Ann. of Math., {\bf 99}(1974), 383-462.

\item{\ref{MG1}} Gregory Margulis,
Discrete groups of motions of manifolds of nonpositive curvature. -
Amer. Math. Soc. Translations, {\bf 109}(1977), 33-45.

\item{\ref{MG2}} \underbar{\phantom{xxxxx}}, Free properly discontinuous
groups of affine transformations.- Dokl. Acad. Sci. USSR, {\bf 272} (1983),
937-940.

%\item{\ref{MM}} Gregory Margulis and George D. Mostow, The differential
%of a quasi-conformal mapping of a Carnot-Carath\'eodory space. - Geom.
%Funct. Anal. {\bf 5} (1995), 402-433.

%\item{\ref{Mas}} Bernard Maskit, Kleinian groups. - Springer-Verlag, 1987.

\item{\ref{Mat}} Rostislav Matveyev, A decomposition of smooth simply-connected
$h$-cobordant 4-mani\-folds. - Preprint, Michigan State Univ. at
E.Lansing, 1995.

\item{\ref{Mil}} John Milnor, On the 3-dimensional Brieskorn manifolds
$M(p,q,r)$.- Knots, groups and \break 3-manifolds, Ann. of Math.
Studies {\bf 84}, Princeton Univ. Press, 1975, 175--225.

\item{\ref{Min}} Robert Miner, Quasiconformal equivalence of spherical CR manifolds.-
Ann. Acad. Sci. Fenn. Ser. A I Math. {\bf 19} (1994), 83-93.

\item{\ref{Mo1}} G.D. Mostow, Strong rigidity of locally symmetric
spaces.- Princeton Univ. Press, 1973.

\item{\ref{Mo2}} \underbar{\phantom{xxxxx}},  On a remarkable class of polyhedra in
complex hyperbolic space. - Pacific J. Math., {\bf 86}(1980), 171-276.

\item{\ref{My}} Robert Myers, Homology cobordisms, link concordances, and
hyperbolic 3-manifolds.- Trans. Amer. Math. Soc.  {\bf 278} (1983), 271-288.

\item{\ref{NR1}} T. Napier and M. Ramachandran, Structure theorems for
complete K\"ahler manifolds and applications to Lefschetz type theorems.-
Geom. Funct. Anal. {\bf 5} (1995), 807-851.

\item{\ref{NR2}} \underbar{\phantom{xxxxx}}, The $L^2$ $\bar\p$-method,
weak Lefschetz theorems, and topology of K\"ahler manifolds.- Preprint,
1996.

\item{\ref{P}} Pierre Pansu, M\'etriques de Carnot-Carath\'eodory et
quasiisom\'etries des espaces symm\'et\-ries de rang un.- Ann. of Math. 
{\bf 129} (1989), 1-60.

\item{\ref{Pr1}} John Parker, Shimizu's lemma for complex hyperbolic space. -
Intern. J. Math. {\bf 3:2} (1992), 291-308.

\item{\ref{Pr2}} \underbar{\phantom{Apanaso}}, Private communication,
1995.
\item{\ref{Ph}} M.B. Phillips, Dirichlet polyhedra for cyclic groups in complex
hyperbolic space. - Proc. AMS, {\bf 115}(1992), 221-228.

\item{\ref{RV}} F.~Raymond and A.T.~Vasquez,
3-manifolds whose universal coverings are Lie groups. -
Topology Appl. {\bf 12} (1981), 161-179.

\item{\ref{Sa}} Nikolai Saveliev, Floer homology and 3-manifold invariants.-
Thesis, Univ. of Oklahoma at Norman, 1995.

\item{\ref{Su1}} Dennis P.~Sullivan, Discrete conformal groups and 
measurable dynamics. - Bull. Amer. Math. Soc. {\bf 6} (1982), 57--73.

\item{\ref{Su2}} \underbar{\phantom{Apanaso}}, Quasiconformal homeomorphisms and dynamics, II: 
Structural stability implies hyperbolicity for Kleinian groups. - Acta Math. 
{\bf 155} (1985), 243--260.

\item{\ref{T1}} Andrew Tetenov, Infinitely generated Kleinian groups in space.
- Siberian Math. J., {\bf 21}(1980), 709-717.

\item{\ref{T2}}  \underbar{\phantom{Apanaso}}, The discontinuity set for
a Kleinian group and topology of its Kleinian manifold. -
Intern. J. Math., {\bf 4:1}(1993), 167-177.

\item{\ref{Th}} William Thurston, The geometry and topology of
three-manifolds. - Lect. Notes, Princeton Univ., 1981.

\item{\ref{To}} Domingo Toledo, Representations of surface groups on complex
hyperbolic space.- J. Diff. Geom. {\bf 29} (1989), 125-133.

\item{\ref{Tu}} Pekka Tukia, On isomorphisms of geometrically finite
Kleinian groups.- Publ. Math. IHES {\bf 61}(1985), 171-214.

\item{\ref{V}} Serguei K. Vodopyanov, Quasiconformal mappings on Carnot
groups.- Russian Dokl. Math. {\bf 347} (1996), 439-442.

\item{\ref{Wa}} C.T.C.~Wall, Geometric structures on compact complex analytic surfaces.
- Topology, {\bf 25} (1986), 119-153.

\item{\ref{Wi}} E.Witten, Monopoles and four-manifolds.- Math. Res.
Lett. {\bf 1} (1994), 769-796.

\item{\ref{Wo}} Joseph A. Wolf, Spaces of constant curvature.- Publ. or Perish,
Berkeley, 1977.

\item{\ref{Ya}} S.T.Yau, Calabi's conjecture and some new results in
algebraic geometry.- Proc. Nat. Acad. Sci {\bf 74} (1977), 1798-1799.

\item{\ref{Yu1}} Chengbo Yue, Dimension and rigidity of quasi-Fuchsian
representations.- Ann. of Math. {\bf 143} (1996)331-355.

\item{\ref{Yu2}}  \underbar{\phantom{Apanaso}}, Mostow rigidity of rank
1 discrete groups with ergodic Bowen-Margulis measure. - Invent. Math.
{\bf 125} (1996), 75-102.

\item{\ref{Yu3}}  \underbar{\phantom{Apanaso}}, Private communication,
Norman/OK, November 1996.

\enddocument